\documentclass{ijuc}

\usepackage[pdftex]{graphics}
\usepackage[utf8x]{inputenc}

\usepackage{tikz}  

\usepackage{amsmath}
\usepackage{fancybox} 
\usepackage{multirow}
\usepackage{amssymb}
\usepackage{permute}  

\newtheorem{thm}{Theorem}[section]
\newtheorem{cor}[thm]{Corollary}

\newtheorem{lemma}[thm]{Lemma}

\newtheorem{conj}[thm]{Conjecture}

\newtheorem{defn}[thm]{Definition}
\newtheorem{example}[thm]{Example}

\newcommand{\Z}{\mathbb Z}
\newcommand{\N}{\mathbb N}

\newcommand{\toffalg}{Toffoli algebra}
\newcommand{\toffalgs}{Toffoli algebras}

\begin{document}

\title{Closed Systems of Invertible Maps}

\author{Tim Boykett\inst{1}
\email{tim.boykett@jku.at and tim@timesup.org}
\footnote{Research supported by the Austrian Science Fund (FWF) research grants P24077, P24285 and Y464-N18.
}
}

\institute{Institute for Algebra, Johannes Kepler University Linz, Austria\\
and
Time's Up Research, Linz, Austria
}

\def\received{}

\maketitle

\begin{abstract}
We  generalise clones, which are sets of functions $f:A^n \rightarrow A$, to sets of maps $f:A^n \rightarrow A^m$. 
We formalise this and develop language that we can use to speak about such maps.
In particular we   look at bijective mappings, which model the logical gates of reversible computation. 
Reversible computation is important for physical (e.g. quantum computation) 
as well as engineering (e.g. heat dissipation) reasons.
We generalise Toffoli's seminal work on reversible computation to multiple valued logics. 
In particular, we show that some restrictions Toffoli found for reversible computation on
alphabets of order 2 do not  apply for odd order alphabets.
For $A$ odd, we can create all invertible mappings from the Toffoli 1- and 2-gates, demonstrating that we can
realise all reversible mappings from four generators.
We  discuss various forms of closure, corresponding to various systems of permitted manipulations.
This leads, amongst other things, to discussions about ancilla bits in quantum computation.

\end{abstract}


\section{Introduction}

The goal of this paper is to generalise Toffoli's fundamental results about reversible computation.
There are several motivations for an interest in reversible computation, some related to physics and engineering, others of a more algebraic nature.
The engineering perspective notes that destroying information creates entropy and thus heat, which is bad for circuitry and wasteful\cite{engineering}.
The physics perspective notes that fundamental physical processes such as inelastic collisions and quantum processes are  reversible\cite{physics}.
An algebraic perspective notes that we can say more about groups than about semigroups owing to the existence of inverses.

 Toffoli  \cite{toff80,tofflncs} introduced what later became known as the Toffoli gate, a simple yet universal basic element for
computation with reversible binary logic.
This paper is mostly concerned with the generalisation of this logic to multiple-valued logics. 
One main result is the generalisation of Toffoli gates to  larger alphabets
and the demonstration that the multiple-valued logics of odd arity are
in some sense more powerful (Section \ref{secgeneration}).

We use language based upon the function algebras of \cite{pk79}, also known as clones\cite{szendrei}.
Mal'cev \cite{malcev} introduced the concept of an iterative algebra, a generalised clone. 
If an iterative algebra includes the projections $\pi_i^n : (x_1,\ldots,x_n) \mapsto x_i$, then it is a \emph{clone}.
In the next sections we will introduce some similar terminology to deal with arbitrary mappings, in particular bijections on $A^n$.
Later we will see that \emph{linear term algebras}, 
the ``reduct of iterative algebras'' that do not have the $\Delta$ operator, play an important role.

In Section \ref{secterm} we introduce mappings and operations upon them, to show that there are several ways to talk about compositional closure.
In particular we show that the algebra of maps can be expressed using a finite signature of finitary operations (Theorem \ref{thmRevcloneEquiv}).
We look at the ways that a given set of maps can generate a closed system of mappings, then at ways of 
embedding other mappings into these closed systems. 
In Section \ref{secrealisation} we talk about \emph{realising} one mapping in a multiclone generated by some other mappings.

We  show that all invertible mappings on a set of odd order can be
generated by four small Toffoli gates and using a definition of generation that is important from
an engineering perspective, show that these are enough for all alphabets.
In Section \ref{secclosures}  look at the various forms of functional closure that arise in the discussion.

\section{Function Terminology}
\label{secterm}

We use $\N$ for the set of positive integers.
Given a set $A$, ${\mathcal O}(A) = \{f:A^n\rightarrow A\mid n \in \N\}$ is the set of all functions on $A$.
We will use the notation of \cite{lehtonen}.
The operations of an iterative algebra are $\tau,\zeta,\Delta,\nabla,*$, defined as follows. Let $f$ be an $n$-ary function, $g$ an $m$-ary function.
The operations $\tau$ and $\zeta$ permute variables, $(\tau f)(x_1,\ldots,x_n) = f(x_2,x_1,x_3,\ldots,x_n)$ 
and $(\zeta f)(x_1,\ldots,x_n) = f(x_2,x_3,\ldots,x_n,x_1)$, with both being the identity on unary functions.
The operation $\Delta$ identifies variables, so $(\Delta f)(x_1,\ldots,x_{n-1}) = f(x_1,x_1,x_2,\ldots,x_{n-1})$, 
while $\Delta$ is the identity on unary functions.
The operation  $\nabla$ introduces a dummy variable, $(\nabla f)(x_1,x_2,\ldots,x_{n+1}) = f(x_2,\ldots,x_{n+1})$.
Lastly we have composition, so
\[(f*g)(x_1,\ldots,x_{n+m-1}) = f(g(x_1,\ldots,x_m),x_{m+1},\ldots,x_{n+m-1}).\]

The \emph{full iterative algebra} is $({\mathcal O}(A); \tau,\zeta,\Delta,\nabla,*)$
and an iterative algebra on $A$ is a subuniverse of this.
If an iterative algebra includes the projections $\pi_i^n : (x_1,\ldots,x_n) \mapsto x_i$, then it is a \emph{clone} on $A$.
The \emph{full linear term algebra} is $({\mathcal O}(A); \tau,\zeta,\nabla,*)$
and a linear term algebra on $A$ is a subuniverse of this.

We use $S_A$ to denote the permutations of a set $A$, $S_n$ to denote the set of bijections on $\{1,\ldots,n\}$ and $S_{(\N)} = \cup_{n\in\N} S_n$.
Juxtaposition multiplies permutations from left to right, so $(1 2 3) (1 2) = (1)(2 3)$, thus in order to 
avoid difficulties when we write permutations as functions, we introduce clarifying parentheses.
Thus if $\alpha = (1 2 3)$ and $\beta = (1 2)$ then $(\alpha \beta)(1) = 1^{\alpha \beta}= (1^\alpha)^\beta= 1$ 
but $\alpha \circ \beta (1) = \alpha(\beta(1)) = (1^\beta)^\alpha=  3$.

Let $A$ be a nonempty set. For $n,m \in \N$, an \emph{$(n,m)$-map} is a mapping $f:A^n\rightarrow A^m$. 
We call $n$ the \emph{arity}, $m$ the \emph{co-arity}.
We write $M_{n,m}(A) = \{f: A^n\rightarrow A^m\}$. 
Then $M(A) = \bigcup_{n,m} M_{n,m}(A)$ is the set of all maps on $A$.
Let $f\in M_{n,m}(A)$, $1\leq i\leq m$. We define $f_i:A^n\rightarrow A$ to be the $i$-th \emph{projection} of $f$,
so we can write $f=(f_1,\dots,f_m)$.
A \emph{function} is an $(n,1)$-map. An $(n,n)$-map, with arity equal to co-arity, will be called \emph{balanced}.
Furthermore, we define $B_{n,m}(A)$ as the set of all $(n,m)$-maps that are bijections, so for $A$ finite $n\neq m$ implies
that $B_{n,m}(A)=\emptyset$ and the bijective maps are balanced.
We will write $B_n(A)$ for the balanced bijections of arity $n$.
Similarly we write $B(A) = \bigcup_{n,m} B_{n,m}(A)$ is the set of all bijections on $A$, $B(A)=\bigcup_{n} B_{n}(A)$ if $A$ is finite.

For example the map $f(x,y) = (x+y,x-y)$ on an abelian group  is a balanced $(2,2)$-map,  it is a bijection if the group is of odd order.

%

%

\begin{defn}
 
 Let $A$ be a set.
 Let $f:A^ m \rightarrow A^n$ ,  $t\leq n$,
 $I=(i_1,\dots,i_t) \in\{1,\ldots,n\}^t$ be a tuple without repetitions. 
 Define $\mu(I,f) \in M_{m,t}(A)$ componentwise as
 $(\mu(I,f))_j = f_{i_j}$, for  all $j \in \{1,\dots,t\}$.
 That is, $\mu(I,f)$ is obtained from $f$ by projecting the output to the components indexed by the tuple $I$.
 
 Let $i\in \{1,\dots,m\}$ and $a\in A$, define $\kappa(i,a,f) \in M_{m-1,n}(A)$ by
 \begin{align*}
  \kappa(i,a,f)(x_1,\dots,x_{m-1}) = f(x_1,\ldots,x_{i-1},a,x_i,\ldots,x_{m-1}).
 \end{align*}
 for all $x_1,\dots,x_{m-1} \in A$. That is, $\kappa(i,a,f)$ is obtained from $f$ by
 substituting the constant $a$ for the $i$-th variable.
 
 We extend $\kappa$ to allow the substitution of constants for several variables simultaneously.
 Let $r\leq m$, $I \in \{1,\dots,m\}^r$  increasing 
 so that $I=(i_1,\dots,i_r)$ with $i_1<i_2<\dots < i_r$ and $a=(a_1,\dots, a_r) \in A^r$, 
define $\kappa(I,a,f) \in M_{m-r,n}(A)$ by the following recursion.
Let $\kappa((),(),f)=f$, for nonempty $I$ define
\begin{align*}
 \kappa((i_1,\dots,i_r),(a_1,\dots,a_r),f) = \kappa(i_1,a_1,\kappa((i_2,\dots,i_{r}),(a_2,\dots,a_{r}),f))
\end{align*}


 Finally we introduce the following shorthand notation.
 For $I \in \{1,\ldots,n\}^t$ and $I^\prime \in \{1,\ldots,m\}^r$, tuples as above,
   and for $o\in A$, define
 \begin{align*}
  f_{I^\prime,I}^o = \mu(I,\kappa(I^\prime,(o,\dots,o),f))
 \end{align*}
Then $f_{I^\prime,I}^o \in M_{m-r,t}(A)$
\end{defn}

For example, take the map $f(x,y) = (x+y,x-y)$ on $\Z_7$, 
then $f_{(2),(1)}^5 (z) = \mu((1),\kappa((2),(5),f))(z) = z+5$ is a unary function on $\Z_7$.

\subsection{Operations}
In this section we look at the operations that allow us to combine  mappings. 
Our main result is that all these combinations can be written as a finite number of finitary operations.
 
\begin{defn}
 \label{defops}
  Let $A$ be a set.
Let $f \in M_{n,s}(A)$, $g\in M_{m,t}(A)$ be  mappings.

Define $i_1 \in {B}_1(A)$ with $i_1(x)=x$, write $i_n \in B_n(A)$ for the identity function on $A^n$.

We  place two mappings next to one another to form a wider mapping.
\begin{align*}
 f\oplus g: A^{n+m} \rightarrow &A^{s+t}\\
 (x_1,\ldots x_{n+m}) \mapsto (&f_1(x_1,\ldots,x_n),\ldots,f_s(x_1,\ldots,x_n),\\
               &g_1(x_{n+1},\ldots,x_{n+m}),\ldots,g_t(x_{n+1},\ldots,x_{n+m}))
\end{align*}

We compose mappings.
For $k \in \N$ with $k\leq n$, $k\leq t$,
\begin{align*}
\circ_k(f,g):A^{n+m-k} \rightarrow   & A^{s+t-k} \\
 (x_1,\ldots,x_{n+m-k}) \mapsto 
   &(f_1(g_1(x_{1},\ldots,x_{m}),\ldots,g_k(x_{1},\ldots,x_{m}),x_{m+1},\ldots,x_{m+n-k}),\\
   & \hspace{12mm}\vdots \\ 
   & \hspace{1mm}f_s(g_1(x_{1},\ldots,x_{m}),\ldots,g_k(x_{1},\ldots,x_{m}),x_{m+1},\ldots,x_{m+n-k}),\\
   & \hspace{1mm}g_{k+1}(x_{1},\ldots,x_{m}), \ldots,g_t(x_{1},\ldots,x_{m}))
\end{align*}
For $k\geq 2$, $\circ_k$ is a partial operation.

We  identify variables.
\begin{align*}
 &\Delta f : A^{n-1} \rightarrow A^s\\
 &\Delta f(x_1,\ldots x_{n-1}) = f(x_1,x_1,\ldots,x_{n-1})
\end{align*}
and define $\Delta f = f $ if $f$ has arity 1.

We introduce dummy variables.
\begin {align*}
 &\nabla f  : A^{n+1} \rightarrow A^s \\
 &\nabla f (x_1,\dots,x_{n+1}) = f(x_2,\dots,x_{n+1})
\end {align*}

The following operations allow us to permute the inputs and outputs of a given mapping.
\begin{align*}
 \tau f (x_1,\ldots,x_n) = &f(x_2,x_1,x_3,\ldots,x_n)\\
 \zeta f (x_1,\ldots,x_n) = &f(x_2,x_3,\ldots,x_n,x_1)\\
 \bar\tau f (x_1,\ldots,x_n) = 
    (&f_2(x_1,x_2,x_3,\ldots,x_n),\\ 
     &f_1(x_1,x_2,x_3,\ldots,x_n),\\
     &\hspace{8mm}\vdots\\
     &f_s(x_1,x_2,x_3,\ldots,x_n))\\
 \bar\zeta f (x_1,\ldots,x_n) = 
    (&f_2(x_1,x_2,x_3,\ldots,x_n),\\ 
     &f_3(x_1,x_2,x_3,\ldots,x_n),\\
     &\hspace{8mm}\vdots\\
     &f_s(x_1,x_2,x_3,\ldots,x_n),\\
     &f_1(x_1,x_2,x_3,\ldots,x_n)).
\end{align*}

In the case that $f$ has arity 1, $\tau f = \zeta f = f$.
In the case that $f$ has co-arity 1,  $\bar\tau f = \bar\zeta f = f$.
\end{defn}


For example, in the case that $n=t$, $f \circ_n g$ is the  composition that we would expect:
\begin{align*}
 f \circ_n g = (&f_1(g_1(x_1,\ldots,x_m),\ldots,g_n(x_1,\ldots,x_m)),\\
 &\hspace{10mm}\vdots\\
 &f_s(g_1(x_1,\ldots,x_m),\ldots,g_n(x_1,\ldots,x_m)))
\end{align*}

Note that there are many relations amongst these operations.
For instance let $f$ have co-arity $n$, then $\bar\tau f = (\tau i_n) \circ_n f$ and $\bar\zeta f = (\zeta i_n) \circ_n f$.
The $\mu$ operator introduced above allows us to say $\mu((2,\dots,n+1),\nabla f) = f$.
We will not examine the axiomatics of these operations in detail.
In order to keep our notation clear, it is important to prove the following result.

\begin{lemma}
\label{lemmacircassoc}
 The operations $\circ_k$ are associative as partial operations, that is, for all $f,g,h \in M(A)$, for all $k,l \in \N$,
 \begin{align}
  f \circ_k (g \circ_l h) = (f \circ_k g) \circ_l h  \label{eqcircassoc}
 \end{align}
 when both sides exist.
\end{lemma}
Proof:
In order to shorten notation, we will write $a_f$ for the arity of $f$ and $c_f$ for the co-arity of $f$.

Suppose both sides of equation $(\ref{eqcircassoc})$ exist.
Then the two sides of  $(\ref{eqcircassoc})$ have equal arity and co-arity, that is
arity $a_f+a_g+a_h - (k+l)$ and co-arity $c_f+c_g+c_h - (k+l)$.
To exist we must have $k \leq \min(a_f,c_g+c_h-l,c_g)$ and $l \leq \min(a_g,c_h,a_f+a_g-k)$.
Then $l\leq c_h$ so $c_g\leq c_g+c_h-l$ so $k \leq \min(a_f,c_g)$ and similarly $l \leq \min(a_g,c_h)$.
We  calculate, for any $1\leq i \leq c_f+c_g+c_h - (k+l)$.
\begin{align*}
 &(f \circ_k (g \circ_l h))_i(x_1,\dots,x_{a_f+a_g+a_h - (k+l)})\\
  &= \begin{cases}
     f_i((g \circ_l h)_1(x_1,\dots,x_{a_g+a_h-l}),\dots,(g \circ_l h)_k(x_1,\dots,x_{a_g+a_h-l}),\nonumber\\
     \hspace{18mm}x_{a_g+a_h-l+1},\dots,x_{a_f+a_g+a_h-(l+k)}) \hspace{8mm} \mbox{ if } i \leq c_f\\
     (g \circ_l h)_{k+i-c_f}(x_1,\dots,x_{a_g+a_h-l}) \hspace{28mm}  \mbox{ if } c_f < i 
    \end{cases}\\
  &= \begin{cases}
     f_i(g_1(h_1(x_1,\dots,x_{a_h}),\dots,h_l(x_1,\dots,x_{a_h}),x_{a_h+1},\dots,x_{a_g+a_h-l}),\dots\\
     \hspace{8mm}\dots,g_k(h_1(x_1,\dots,x_{a_h}),\dots,h_l(x_1,\dots,x_{a_h}),x_{a_h+1},\dots,x_{a_g+a_h-l}),\nonumber\\
     \hspace{8mm}x_{a_g+a_h-l+1},\dots,x_{a_f+a_g+a_h-(l+k)})     \hspace{18mm}\mbox{ if }i \leq c_f\\
     g_{k+i-c_f}(h_1(x_1,\dots,x_{a_h}),\dots,h_l(x_1,\dots,x_{a_h}),x_{a_h+1},\dots,x_{a_g+a_h-l}) \nonumber\\
     \hspace{78mm}\mbox{ if }c_f < i \leq c_f+c_g-k \\
     h_{l+i-(c_f+c_g-k)}(x_1,\dots,x_{a_h})   \hspace{37mm} \mbox{ if }c_f+c_g-k < i
    \end{cases}
\end{align*}
while
\begin{align*}
 &((f \circ_k g) \circ_l h)_i(x_1,\dots,x_{a_f+a_g+a_h - (k+l)})\\
 &=\begin{cases}
    (f \circ_k g)_i(h_1(x_1,\dots,x_{a_h}),\dots,h_l(x_1,\dots,x_{a_h}),\nonumber\\
     \hspace{28mm}x_{a_h+1},\dots,x_{a_f+a_g+a_h-(k+l)}) \hspace{2mm}\mbox{ if } i \leq c_f+c_g-k \\
    h_{l+i-(c_f+c_g-k)}(x_1,\dots,x_{a_h})    \hspace{31mm} \mbox{ if }i > c_f+c_g-k
   \end{cases}
\end{align*}

\begin{align*}
 &=\begin{cases}
    f_i(g_1(h_1(x_1,\dots,x_{a_h}),\dots,h_l(x_1,\dots,x_{a_h}),x_{a_h+1},\dots,x_{a_g+a_h-l}),\dots\\
    \hspace{8mm}\dots,g_k(h_1(x_1,\dots,x_{a_h}),\dots,h_l(x_1,\dots,x_{a_h}),x_{a_h+1},\dots,x_{a_g+a_h-l}),\nonumber\\
    \hspace{8mm}x_{a_g+a_h-l+1},\dots,x_{a_f+a_g+a_h-(k+l)})    \hspace{13mm}\mbox{ if } i \leq c_f \\
    g_{k+i-c_f}(h_1(x_1,\dots,x_{a_h}),\dots,h_l(x_1,\dots,x_{a_h}),x_{a_h+1},\dots,x_{a_g+a_h-l}) \nonumber\\
     \hspace{70mm}\mbox{ if } c_f < i \leq c_f+c_g-k \\
    h_{l+i-(c_f+c_g-k)}(x_1,\dots,x_{a_h})     \hspace{28mm}\mbox{ if } i > c_f+c_g-k
   \end{cases}
\end{align*}
These can be seen to be equal and we are done.
\hfill$\Box$

Note that one side of the associativity equation in the previous lemma can exist
without the other existing, for instance if $k=l=2$, $f,h \in M_{2,3}(A)$ and $g \in M_{2,1}(A)$ then
 $f \circ_k (g \circ_l h)$ exists while $ (f \circ_k g) \circ_l h$ does not.

We can look at a general form of mapping composition corresponding to clones.
Remember that to be closed under a 0-ary operation such as $i_1$ means to include that constant.
\begin{defn}
 Let $A$ be a set.
A collection of  mappings on $A$ is called a  \emph{multiclone} on $A$ if it is closed under  
the operations $\{i_1,\oplus,\Delta,\nabla,\tau,\zeta\}$ and the partial operations $\{\circ_i\mid i\in \N\}$.
Let $F\subseteq M(A)$ be a collection of mappings. 
Then we write $C_\Delta(F)$ for the closure under the operations,\ i.e. the \emph{multiclone generated by} $F$.
\end{defn}

However in order to focus on the realm of bijective mappings,
we are predominately interested in closure without the $\Delta$ and $\nabla$ operations.

\begin{defn}
 Let $A$ be a set.
A collection of  mappings on $A$ is called a \emph{\toffalg}  on $A$ if it is closed under  
the operations $\{i_1,\oplus,\tau,\zeta\}$ and the partial operations $\{\circ_i\mid i\in \N\}$.
Let $F\subseteq M(A)$ be a collection of mappings. 
Then we write $C(F)$ for the closure under the operations,\ i.e. the \emph{\toffalg{}  generated by} $F$.
\end{defn}

Note that by the comments following Definition \ref{defops} we can express the $\bar \tau$ and $\bar \zeta$
operations in terms of these operations, thus not losing any expressive power.

We note that there are similarities to the Network Algebra used in e.g. \cite{stefanescu13,stefanescu2000}.
We note that the expression Toffoli Algebra has been used in \cite{jensen14} in a related but distinct form.

These operations reflect what Toffoli \cite{toff80} calls a \emph{combinatorial network}, 
as well as corresponding to the idea of a \emph{read once function} as used in in \cite{coucierolehtonen}.
In particular, every output is used at most once as an input to another mapping, 
avoiding the duplication of variables or so-called fan-out,
where one data signal is spread to two or more outputs. The mapping
$\phi_n \in M_{1,n}(A)$ with $\phi_n(a)=(a,\dots,a)$ is a \emph{fan-out} mapping.
Also, no input or output is thrown away. 
As a connection between these concepts, we have the following.

\begin{lemma}
\label{lemmfan}
 Let $F$ be a \toffalg. Then $F$ is a multiclone iff $\phi_2 \in F$ and $g\in F$ where
 $\phi_2(x)=(x,x)$, $g(x,y)=y$.
\end{lemma}
Proof:
$(\Rightarrow)$:  We know that $i_2=i_1\oplus i_1\in F$. We know $F$ is a multiclone,
so $\Delta i_2 \in F$. But $\Delta i_2(x_2)=(x_2,x_2) = \phi_2(x_2)$ so $\phi_2 \in F$.
Similarly $\nabla i_1 = g$.

$(\Leftarrow)$: 
Suppose $\phi_2\in F$. Let $h\in F\cap M_{n,m}$ with $n \geq 2$.
Then 
\begin{align*}
 (h \circ_2 \phi_2)(x_1,\ldots,x_{n-1}) &= h(x_1,x_1,\ldots,x_{n-1})\\
 &= \Delta h(x_1,\ldots x_{n-1})
\end{align*}
so $F$ is closed under $\Delta$.

Let $h\in F\cap M_{n,m}$.
We calculate that $h \circ_1 g \in M_{n+2-1,m+1-1}=M_{n+1,m}$
\begin{align*}
 (h \circ_1 g)(x_1,\dots,x_{n+1}) &= h(g_1(x_1,x_2),x_3,\dots,x_{n+1}) \\
  & = h(x_2,\dots,x_{n+1}) \\
  & = \nabla h(x_1,\dots,x_{n+1})
\end{align*}
so $F$ is also closed under $\nabla$, thus $F$ is a multiclone.
\hfill$\Box$


\begin{defn}
A \toffalg{}  is called a \emph{reversible \toffalg}  if each of its elements element is bijective.
\end{defn}


Note that $i_1\oplus i_1\oplus \ldots\oplus i_1=i_{n}$.
We will also use $\oplus$ to concatenate tuples, 
$(x_1,\dots,x_n)\oplus(y_1,\dots,y_m) = (x_1,\dots,x_n,y_1,\dots,y_m)$.
Let $\alpha \in S_n$ be a permutation.
Then let $\pi_\alpha$ be the function 
$\pi_\alpha (x_1,\ldots,x_n) = (x_{\alpha^{-1}1},\ldots, x_{\alpha^{-1}n})$.

\begin{example}
 The set $\{\pi_\alpha\mid \alpha \in S_n,\,n\in \N\}$ is a reversible \toffalg, in fact it is the unique  minimal  \toffalg{} on any set $A$.
 \end{example}

Other multiclones and reversible \toffalgs{} arise from linear algebra.

\begin{example}
 Let $A$ be a field. Then the collection of linear mappings on vector spaces over $A$ form a multiclone.
 The collection of invertible linear mappings $\cup_{n\in \N} GL(n,A)$ forms a reversible \toffalg, but not a multiclone. 
 The mappings described by permutation matrices form the smallest Toffoli subalgebra and as a subalgebra of a reversible
 \toffalg{} it is also reversible.
\end{example}

\begin{defn}
Let $f,g$ be  mappings, $f$ have arity $n$, $g$ have co-arity $m$.
 Define $f \bullet g = f \circ_k g$ where $k=\min(m,n)$.
\end{defn}

This is a total operation.

\begin{lemma}
Let $f,g$ be  mappings, $f$ have arity $n$, $g$ have co-arity $m$.
Then
\begin{align*}
f\bullet g =
\begin{cases}
  (f \oplus i_{m-n}) \bullet g =  (f \oplus i_{m-n}) \circ_m g & \mbox{ if }n < m\\
  f  \circ_m g & \mbox{ if } n=m \\
   f  \bullet (g\oplus i_{n-m}) =  f  \circ_n (g\oplus i_{n-m}) & \mbox{ if }n > m
\end{cases}
\end{align*}
Also for any $k\in \N$, $(f\bullet g)\oplus i_k = (f \oplus i_k) \bullet (g \oplus i_k)$.
\end{lemma}
Proof:
Let $f$ have co-arity $s$ and $g$ have arity $t$.

Suppose $n<m$, $f \bullet g = f \circ_n g$ has arity $n+t-n=t$. We see that 
 $f \oplus i_{m-n}$ has arity $n+m-n=m$ so $(f \oplus i_{m-n}) \bullet g =(f \oplus i_{m-n}) \circ_m g$ holds.
Then 
\begin{align*}
 &(f \bullet g) (x_1,\dots,x_{t})\\
  &= (f \circ_n g) (x_1,\dots,x_{t})\\
  & = (f_1(g_1(x_1,\dots,x_{t}),\dots,g_n(x_1,\dots,x_{t})),\dots \nonumber\\
  & \hspace{18mm}\dots,f_s(g_1(x_1,\dots,x_{t}),\dots,g_n(x_1,\dots,x_{t})),\nonumber\\
  & \hspace{18mm}         g_{n+1}(x_1,\dots,x_{t}),\dots,g_{m}(x_1,\dots,x_{t})) \\
  & =(f \oplus i_{m-n})(g_1(x_1,\dots,x_{t}),\dots,g_m(x_1,\dots,x_{t}))\\
  & = ((f \oplus i_{m-n}) \circ_m g) (x_1,\dots,x_{t})
\end{align*}
which is what we wanted.

The second claim is immediate from the definition of $f \bullet g$.

The third claim is proved analogously to the first.

For the fourth claim we start with the case that $n=m$.
We calculate
\begin{align*}
 &((f\bullet g)\oplus i_k) (x_1,\dots,x_{t+k})\\
 &=(f \bullet g)(x_1,\dots,x_t) \oplus (x_{t+1},\dots,x_{t+k}) \\
 &=(f \circ_n g)(x_1,\dots,x_t) \oplus (x_{t+1},\dots,x_{t+k}) \\
 &=(f_1(g_1(x_1,\dots,x_t),\dots,g_n(x_1,\dots,x_t)),\dots\nonumber\\
 & \hspace{8mm}\dots,f_s(g_1(x_1,\dots,x_t),\dots,g_n(x_1,\dots,x_t)),x_{t+1},\dots,x_{t+k}) \\
 &= ((f \oplus i_k) \circ_n g) (x_1,\dots,x_{t+k})\\
 &= ((f \oplus i_k) \circ_{n+k} (g\oplus i_k)) (x_1,\dots,x_{t+k})\\
 &= ((f \oplus i_k) \bullet (g\oplus i_k)) (x_1,\dots,x_{t+k})
\end{align*}
so the claim holds in this case.
Now note that if $n<m$ then $f\oplus i_{m-n}$ has arity $m$ so 
$(f \bullet g) \oplus i_k = ((f\oplus i_{m-n}) \bullet g)\oplus i_k = (f \oplus i_{m-n} \oplus i_k) \bullet (g \oplus i_k) = 
(f \oplus i_k) \bullet (g \oplus i_k)$ so this case is done. Similarly we argue in the case that $n>m$ and thus we are done.
\hfill$\Box$

\begin{lemma}
 The $\bullet$ operation is associative, that is, for all $f,g,h \in M(A)$,
 $f \bullet (g\bullet h) = (f \bullet g) \bullet h$.
\end{lemma}
Proof:
In order to simplify notation, we write $a_f$ for the arity of $f$ and $c_f$ for the co-arity of $f$.
Then $f \bullet (g\bullet h) = f \circ_k (g \circ_l h)$ where $k=\min(a_f,c_g+c_h-l)$ and $l=\min(a_g,c_h)$.
Similarly $(f \bullet g) \bullet h = ( f \circ_{\bar k} g) \circ_{\bar l} h$ where $\bar k = \min(a_f,c_g)$ and $\bar l = \min(a_f+a_g-\bar k, c_h)$.

Suppose $a_f \leq c_g$.
Then $\bar k = k = a_f$ so $\bar l = \min(a_g,c_h)=l$ so by the associativity of $\circ_i$ for all $i$, Lemma \ref{lemmacircassoc}, we are done.

Similarly if $c_h \leq a_g$, $l=\bar l = c_h$ so $k = \min(a_f,c_g)=\bar k$ and we are done.

Suppose that neither of these cases applies, that is,  $a_f> c_g$ and $c_h>a_g$.
Then
\begin{align}
 (f \bullet g)\bullet h = (f \circ_{a_f} (g \oplus i_{a_f-c_g})) \bullet h \label{asseq1}\\
 f \bullet (g \bullet h) = f \bullet ((g \oplus i_{c_h-a_g}) \circ_{c_h} h) \label{asseq2}
\end{align}
Let $a_1$ be the arity of $f \bullet g = f \circ_{a_f} (g \oplus i_{a_f-c_g})$, that is $a_1= a_g+a_f-c_g$.
Let $c_2$ be the co-arity of $g \bullet h =(g \oplus i_{c_h-a_g}) \circ_{c_h} h$, that is $c_2=c_g+c_h-a_g$.
Note that $a_1-c_h=a_f-c_2$.

We proceed by cases.

Case $a_1 < c_h$ so $ c_2 > a_f$. Then
\begin{align*}
 (\ref{asseq1}) &= ((f \circ_{a_f} (g \oplus i_{a_f-c_g}))\oplus i_{c_h-a_1}) \circ_{c_h} h \\
  &= ( (f \oplus i_{c_h-a_1})  \circ_{a_f+c_h-a_1} (g \oplus i_{a_f-c_g}\oplus i_{c_h-a_1}))\circ_{c_h} h\\
  &= ( (f \oplus i_{c_2-a_f})  \circ_{c_2} (g \oplus i_{c_h-a_g}))\circ_{c_h} h\\
  &=  (f \oplus i_{c_2-a_f})  \circ_{c_2}(( g \oplus i_{c_h-a_g})\circ_{c_h} h)\\
  &=  f   \bullet( g \bullet h) = (\ref{asseq2})
\end{align*}

Case $a_1 = c_h $ so $c_2 = a_f$. Then $a_f-c_g=c_h-a_g$ and
\begin{align*}
 (\ref{asseq1}) 
  &= (f \circ_{a_f} (g \oplus i_{a_f-c_g})) \circ_{c_h} h \\
  &=  f   \circ_{c_2}(( g \oplus i_{c_h-a_g})\circ_{c_h} h)\\
  &=  f   \bullet( g \bullet h) = (\ref{asseq2})
\end{align*}

Case $a_1 > c_h$ so $ c_2 < a_f$. Then
\begin{align*}
 (\ref{asseq1}) &= (f \circ_{a_f} (g \oplus i_{a_f-c_g})) \circ_{a_1} (h\oplus i_{a_1-c_h}) \\
  &= f \circ_{a_f} ((g \oplus i_{a_f-c_g}) \circ_{a_1} (h\oplus i_{a_1-c_h})) \\
  &= f \circ_{a_f} ((g \oplus i_{a_f-c_g}) \bullet (h\oplus i_{a_1-c_h})) \\
  &= f \circ_{a_f} (((g \oplus i_{c_h-a_g}) \bullet h)\oplus i_{a_f-c_2}) \\
  &= f \bullet ((g \oplus i_{c_h-a_g}) \bullet h) = (\ref{asseq2})
\end{align*}
and we are done.
\hfill$\Box$

\begin{thm}\label{thmRevcloneEquiv}
 Let $F$ be a collection of  mappings on a set $A$. Then the following are equivalent:
 \begin{enumerate}
  \item $F$ is a \toffalg.
  \item $F$ is closed under $\{i_1,\oplus,\tau,\zeta\} \cup\{ \circ_k\mid k\in \N\}$.
  \item $F$ is closed under $\{\oplus\}\cup\{ \pi_\alpha\mid \alpha \in S_{(\N)}\} \cup\{ \circ_k \mid  k \in \N\}$.
  \item $F$ is closed under $\{i_1,\oplus,\tau,\zeta,\bullet \}$.
 \end{enumerate}
\end{thm}
Proof: The second case is a writing of the definition of the first.
We demonstrate that the others are equivalent by showing that each collection of operations can be built from the others as terms.

$3\Rightarrow 2:$ 
In this signature, $i_1$ is $\pi_\alpha$ for $\alpha$ the identity permutation on $\{1\}$.
$\oplus$ is the same, as are $\circ_k$.
Let $f \in M_{m,n}$.
We see that $\tau f  = f \circ_m \pi_{(1,2)}$,  the permutation $(1,2)$ acting upon $\{1,\ldots,m\}$,
and $\zeta f = f \circ_m \pi_{(m,\ldots,2,1)}$,  the permutation $(m,\ldots,2,1)$ acting upon $\{1,\ldots,m\}$.

$2\Rightarrow 4:$ 
From the definition of $\bullet$, we know that it can be written in terms of $\circ_k$ for any given arguments.

$4\Rightarrow 3:$
Let $\alpha$ be a permutation on $\{1,\ldots,n\}$.
Let $\alpha = \alpha_1\alpha_2\ldots\alpha_k$ for $\alpha_i \in \{(1,2),(1,2,\ldots,n)\}$. 
Let $\beta_i = \tau$ if $\alpha_i=(1,2)$,
otherwise $\beta_i=\zeta$.
Then $\pi_\alpha=\beta_k\ldots\beta_1 (i_1\oplus \ldots \oplus i_1)$ with $i_1$ repeated $n$ times.

Let $f,g$ be  mappings, $f\in M_{n,s}(A)$, $g\in M_{m,t}(A)$, $k\leq t$, $k\leq n$.
Then we claim that
\begin{align}
 f\circ_k g = (f \oplus i_{t-k} ) \bullet \pi_\alpha \bullet (g \oplus i_{n-k}) \label{eqfullcomp}
\end{align}
 with 
\[\alpha = \pmtv[1\ldots {k}{k+1}\ldots {t}{t+1}\ldots{t+n-k}] {{k+1}{n+1} {t}{t+n-k} {t+1}{k+1} {t+n-k}{n} } \in S_{n+t-k}.\]

The left hand side of $(\ref{eqfullcomp})$ has arity $n+m-k$, co-arity $s+t-k$ and is defined by
\begin{align*}
 f \circ_k g (x_1,\dots, x_{n+m-k}) = 
 \begin{cases}
  f_i(g_1(x_1,\dots,x_m),\dots,g_k(x_1,\dots,x_m),\\
    \hspace{18mm}x_{m+1},\dots,x_{m+n-k}) & \mbox{ if }i \leq s \\
  g_{k+i-s}(x_1,\dots,x_m) &\mbox{ if } s<i
 \end{cases}
\end{align*}

On the right hand side of $(\ref{eqfullcomp})$, 
we note that the co-arity of $g \oplus i_{n-k}$ and the arity of $f \oplus i_{t-k}$ are both
$t+n-k$, the same as the arity and co-arity of $\pi_\alpha$.
Because of the structure of the permutation $\alpha$,
\begin{align*}
 (\pi_\alpha \bullet (g \oplus i_{n-k}))_j (x_1,\dots,x_{m+n-k}) = 
 \begin{cases}
  g_j(x_1,\ldots,x_m) & \mbox{ if }j \leq k\\
  x_{m-k+j} & \mbox{ if }k < j \leq n \\
  g_{k+j-n}(x_1,\ldots,x_m) &\mbox{ if }n < j \leq t+n-k
 \end{cases}
\end{align*}
We can then calculate the right hand side of $(\ref{eqfullcomp})$. For $j\leq s$, the $j$th entry is:
\begin{align*}
 &f_j\bullet \pi_\alpha \bullet (g \oplus i_{n-k})(x_1,\ldots,x_{n+m-k}) = \\
     & \hspace{15mm}f_j(g_1(x_1,\ldots,x_m),\ldots,g_k(x_1,\ldots,x_m),x_{m+1},\ldots,x_{m+n-k})
\end{align*}
while for $s<j\leq s+t-k$, the $j$th entry is the output of the identity map, 
so noting that the preimage of $n+j-s$ by $\alpha$ is $k+j-s$:
\begin{align*}
&((f \oplus i_{t-k} )\bullet \pi_\alpha \bullet (g \oplus i_{n-k}))_j(x_1,\ldots,x_{n+m-k}) = \\
 &\hspace{10mm}( \pi_\alpha \bullet (g \oplus i_{n-k}))_{n+j-s}(x_1,\ldots,x_{n+m-k}) = g_{k+j-s}(x_1,\ldots,x_m)
\end{align*}
By comparing entries, we see that the left hand side and the right hand side agree at all entries and are thus equal.
\hfill$\Box$

The following result enables us to  see that the operations of a \toffalg{}  do not destroy reversibility,
so we know that $B(A)$ is the largest reversible \toffalg{}  on $A$.

\begin{lemma}
\label{lemmaBclosed}
 Let $A$ be a set.
 Let $F\subseteq B(A)$. Then $C(F) \subseteq B(A)$.
\end{lemma}
Proof:
We only need to check that $B(A)$ is a \toffalg.
The operations $i_1,\tau,\zeta$ are invertible.
Let $f,g\in B(A)$.
The inverse of $f\oplus g $ is $f^{-1} \oplus g^{-1}$.
If the co-arity of $g$ is the same as the arity of $f$, then $(f \bullet g)^{-1} = g^{-1} \bullet f^{-1}$.
If the co-arity of $g$ is less than the arity of $f$, then note that $f\bullet g = f \bullet (g\oplus i_s)$ where $s$
is the difference between the arity of $f$ and the co-arity of $g$. Then the co-arity of $g\oplus i_s$ equals the arity of $f$
and we are done. Similarly if the arity of $f$ is less that the co-arity of $g$.

Thus all of the operations of a \toffalg{}  map bijections to bijections, so we are done.
\hfill$\Box$

We  know that there is a finite signature for \toffalgs, 
 so we can apply the standard tools and techniques of universal algebra,
 such as the following.

 \begin{cor}
 \label{corAlgLattice}
  Let $A$ be a set. Then the set of \toffalgs{} on $A$, ordered by inclusion, is an algebraic lattice.
  The set of reversible \toffalgs{} on $A$, ordered by inclusion, is an algebraic lattice.
 \end{cor}
Proof: The set of \toffalgs{} on $A$ is the set of subuniverses of the algebra $(M(A); i_1,\oplus,\tau,\zeta,\bullet)$.
These operations are all finitary, so
by \cite{BS}[Cor 3.3] we have our result.

Similarly reversible \toffalgs{} are subuniverses of $(B(A); i_1,\oplus,\tau,\zeta,\bullet)$ and we are done.
\hfill$\Box$

We note also that the definition of a multiclone adds only two unary operations to the definition of a \toffalg, so we obtain the following
as a corollary to Theorem \ref{thmRevcloneEquiv} and Corollary \ref{corAlgLattice}.
\begin{cor}
  Let $F$ be a collection of  mappings on a set $A$. Then the following are equivalent:
 \begin{enumerate}
  \item $F$ is a multiclone.
  \item $F$ is closed under $\{i_1,\oplus,\tau,\zeta,\Delta,\nabla\} \cup\{ \circ_k\mid k\in \N\}$.
  \item $F$ is closed under $\{\oplus,\Delta,\nabla\}\cup\{ \pi_\alpha\mid \alpha \in S_{(\N)}\} \cup\{ \circ_k \mid  k \in \N\}$.
  \item $F$ is closed under $\{i_1,\oplus,\tau,\zeta,\Delta,\nabla,\bullet \}$.
 \end{enumerate}
 The set of multiclones, ordered by inclusion, is an algebraic lattice.
\end{cor}

\section{Realisation}
\label{secrealisation}

In this section we look at the ways in which one set of mappings can be found within another.

\begin{thm}
\label{thmToffFund}
Let $A$ be a finite set, $g:A^m \rightarrow A^n$, $o\in A$ and
 let $r = \max(m,n+\lceil\log_{\vert A\vert} max_{a\in A^n} \vert g^ {-1} (a)\vert \rceil)$.
 There exists an invertible $f \in B_r(A)$,
  such that $g = f_{\theta_1,\theta_2}^o$
  with $\theta_1=(m+1,\ldots,r), \theta_2 = (1,\ldots,n)$.
 Furthermore $\max(m,n) \leq r \leq n+m$ and these bounds are achieved.
\end{thm}
Proof:
Select $o\in A$ arbitrary but fixed.

Suppose that $r=n$. 
This occurs iff $\lceil\log_{\vert A\vert} max_{a\in A^n} \vert g^ {-1} (a)\vert \rceil =0$
and $m \leq n$. 
Now $\lceil\log_{\vert A\vert} max_{a\in A^n} \vert g^ {-1} (a)\vert \rceil =0$ iff 
$max_{a\in A^n} \vert g^ {-1} (a)\vert=1$ iff
$g$ is
injective. If $m = n$ then $g$ is already bijective and we are done.
If $m < n$ then define the partial map $f$ on $A^r$
by $f(x_1,\dots,x_m,o,\dots,o) = g(x_1,\dots,x_m)$ for all $x_1,\dots,x_m \in A$.
Then
$f$  is a partial injective map on $A^{r}$, thus it can be completed to a bijection $f$ of $A^{r}$.
Then $f_{\theta_1,\theta_2}^o = \kappa((m+1,\dots,n),o,f)=g$
and we are done.

Now suppose that $r >n$.
For each $a \in A^n$, let $x^{(1)},\ldots,x^{(k)} \in A^m$ be a lexicographical ordering of $g^{-1}(a)$.
Let $b^{(1)},\ldots,b^{(k)}$ be the first $k$ elements of $A^{r-n}$ in lexicographical order.
Define $f(x^{(i)}_{1},\ldots,x^{(i)}_{m},o,\ldots,o) = (a_1,\ldots,a_n,b^{(i)}_{1},\ldots,b^{(i)}_{(r-n)})$ for each $i \in \{1,\ldots,k\}$.

Now $f$ is a partial injective map on $A^{r}$, thus it can be completed to a bijection of $A^{r}$.

For any $z \in A^m$, let $b=g(z_1,\ldots,z_m)$.
Then $(z_1,\ldots,z_m) \in g^{-1}(b)$, so
$f(z_1,\ldots,z_m,o,\ldots,o) = (b_1,\ldots,b_n,c_1,\ldots,c_m)$ for some $c_1,\ldots,c_m\in A$.

Thus
$f^o_{\theta_1,\theta_2}(z) =\mu(\theta_2,\kappa(\theta_1,(o,\ldots,o),f)(z) = b = g(z)$,
showing $f_{\theta_1,\theta_2}^o = g$.

The average value of $\vert g^{-1}(a)\vert$ is $\frac{\vert A\vert^m}{\vert A\vert^n} = \vert A\vert^{m-n}$.
Thus we know that $\lceil\log_{\vert A\vert} max_{a\in A^n} \vert g^ {-1} (a)\vert\rceil \geq (m-n)$.
In the case that $m \geq n$, this value can be minimised if all $\vert g^{-1}(a)\vert$ are equal to the average value, 
so $\vert g^{-1}(a)\vert = \vert A\vert^{m-n}$, so the minimum value of $r$ is $m=\max(m,n)$.
In the case that $m < n$, this value is minimised if $g$ is injective, so $\vert g^{-1}(a)\vert \leq 1$ 
 and thus  $\lceil\log_{\vert A\vert} max_{a\in A^n} \vert g^ {-1} (a)\vert\rceil =0$
so $r=\max(m,n+0)$.
The maximum value is obtained if $g$ is constant, so there is some $a\in A^n$ such that $g^{-1}(a) = A^m$.
In this case, 
$\log_{\vert A\vert} max_{a\in A^n} \vert g^ {-1} (a)\vert = m$ and thus $r=n+m$.
\hfill$\Box$

\begin{cor}[ \cite{toff80} Theorem 4.1]
 \label{corToffFund}
Let $A$ be a finite set.
 For all $g:A^m \rightarrow A^n$ there exists $r  \leq n$, invertible $f:A^{m+r} \rightarrow A^{m+r}$, $o\in A$
 and $\theta_1 \in \{1,\ldots,m+r\}^m$, $\theta_2 \in \{1,\ldots,m+r\}^n$ such that $g = f_{\theta_1,\theta_2}^o$.
\end{cor}
Proof:
Let $r_0$ be the value of $r$ obtained in the theorem above.
There are two possibilities for $r_0$. If $r=m$ then we take $r=0<n$ in this theorem and we are done.
If $r_0=n+\lceil\log_{\vert A\vert} max_{a\in A^n} \vert g^ {-1} (a)\vert \rceil$ then
note that $g^ {-1} (a) \subseteq A^m$ so $\vert g^ {-1} (a)\vert \leq \vert A\vert ^m $.
So $r_0 \leq n+m$ so $r \leq n$ in this theorem and we are done.
\hfill$\Box$

\subsection{Forms of Realisation}

We now look at several ideas about generation, corresponding to different forms of closure.
While the simplest definition talks about the multiclone generated from a set of mappings, 
Toffoli introduces some interesting ideas about more general, engineering inspired forms of realisation.
The following is a translation of his  definitions to the
language we have developed here.

\begin{defn}
 Let $A$ be a set, $F\subseteq M(A)$. 
 Then $K(F)$ is the set of all functions we can obtain from $F$ by inserting constants using the $\kappa$ operation.
Hence $K(F)$ the closure under all possible applications  of $\kappa$, including none.
 
 $S(F)$ is the set of all submappings of mappings in $F$.
Then $S(F)$ is the closure under all possible applications of $\mu$. 

 $R(F)$ is the collection of bijections in $F$, i.e. $R(F) = F \cap B(A)$.
\end{defn}

\begin{defn}
Let $F$ be a collection of mappings on a set $A$. 
A mapping $g\in M(A)$ is said to be \emph{isomorphically realised} by $F$\cite[p. 8]{toff80}, 
or \emph{realized without ancilla bits} \cite{yangetal}, if $g\in C(F)$.

 Let $A$ be a set, $F \subseteq M(A)$ and $g \in M_{m,n}(A)$. 
 \begin{itemize}
\item  We say $g$ is \emph{realised by} $F$ if there exists some $f\in C(F)\cap M_{l,k}(A)$ and $a_1,\ldots,a_{l-m}\in A$ such that for all $i \in \{1,\ldots,n\}$ and for all $x_1,\ldots,x_m\in A$:
 \begin{align}
  g_i(x_1,\ldots,x_m) = f_i(x_1,\ldots,x_m,a_1,\ldots,a_{l-m}) \label{eqrealize}
 \end{align}
 Equivalently, we can say $g \in SKC(F)$.
\item  We say $g$ is \emph{realised with no garbage} if  there exists some $f \in M_{l,n}(A)\cap C(F)$ satisfying $(\ref{eqrealize})$, 
equivalently $g\in KC(F)$. 
Otherwise $g$ is \emph{realised with garbage}.
\item We say  $g$ is \emph{realised with no constants}  if there exists some $f \in M_{m,k}(A)\cap C(F)$  satisfying  $(\ref{eqrealize})$, 
equivalently $g\in SC(F)$.  
Otherwise $g$ is \emph{realised with constants}.
\end{itemize}
\end{defn}

A mapping  is isomorphically realised iff it is realised with no constants and no garbage, that is,
that there exists a single function that satisfies both claims.
This should not be confused with the claim that for any $F \subseteq M(A)$, $C(F) = SC(F) \cap KC(F)$, which
remains open. The inclusion $\subseteq$ is clear since $C(F) \subseteq SC(F)$ and $C(F) \subseteq KC(F)$ for any $F\subseteq M(A)$,
the other inlusion is unclear.

However  we obtain the following.
\begin{lemma}
 Let $A$ be a  set, $F \subseteq M(A)$ a collection of balanced maps, then $C(F) = SC(F) \cap KC(F)$.
 In particular, let $A$ be a finite set, $F \subseteq B(A)$, then $C(F) = SC(F) \cap KC(F)$.
\end{lemma}
Proof:
By the comment above, we only need show that $SC(F) \cap KC(F) \subseteq C(F)$.
Suppose $g\in SC(F) \cap KC(F)$ and $g \in M_{m,n}(A)$.
Also note that all the operations of a \toffalg{}  map balanced maps to balanced maps, so $C(F)$ consists only of balanced maps.

Since $g\in KC(F)$ we know there exists some $h\in C(F)\cap M_{l,n}(A)$ and $a_1,\dots,a_{l-m} \in A$ such that
for all $i \in \{1,\dots,n\}$, $g_i(x_1,\dots,x_m) = h_i(x_1,\dots,x_m,a_1,\dots,a_{l-m})$.
Because $h$ is balanced, $l=n$, so $m \leq n$.

Since $g\in SC(F)$ we know there exists some $\bar h\in C(F)\cap M_{m,k}(A)$  such that
for all $i \in \{1,\dots,n\}$, $g_i(x_1,\dots,x_m) = \bar h_i(x_1,\dots,x_m)$.
Because  $\bar h$ is balanced, $m=k$, so $n \leq k$.

Combining these we see that $m=n=k=l$ so $g=h=\bar h$ and $g\in C(F)$ so we are done.

We know that when $A$ is finite, bijections in $B(A)$ are balanced, so we obtain the second result.
\hfill$\Box$

This language allows us to rephrase Corollary \ref{corToffFund} as:
\begin{cor}
 Let $A$ be a finite set. 
 For all $g\in M_{m,n}(A)$, $g$ is realised by $B_{m+n}(A)$.
\end{cor}

There is a special class of realisations that have constants and garbage, 
but in some way do not use them up, the garbage representing the constants.

\begin{defn}
 Let $A$ be a set, $F \subseteq M(A)$ and $g \in M_{m,n}(A)$. 
 Then $g$ is \emph{realised with temporary storage by} $F$ if there exists some $f\in C(F)$, 
 $f\in M_{l,k}(A)$ and $a_1,\ldots,a_{l-m}\in A$ such that $n+l-m=k$ and
 \begin{align}
  \forall i \in \{1,\ldots,n\}:&\, g_i(x_1,\ldots,x_m) = f_i(x_1,\ldots,x_m,a_1,\ldots,a_{l-m})\label{eqt1}\\
  \forall i \in \{1,\ldots,l-m\}&\, \forall x_1,\ldots,x_m \in A:\nonumber\\ 
      &\, f_{n+i}(x_1,\ldots,x_m,a_1,\ldots,a_{l-m}) = a_i  \label{eqt2}
 \end{align}
\end{defn}

\begin{defn}
 Let $A$ be a set, $F \subseteq M(A)$ and $g \in M_{m,n}(A)$. 
 Then $g$ is \emph{realised with strong temporary storage by} $F$ if there exists some $f\in C(F)$, $f\in M_{l,k}(A)$ 
 and $a_1,\ldots,a_{l-m}\in A$ such that $n+l-m=k$ and
 \begin{align*}
  \forall i \in \{1,\ldots,n\}:&\,  g_i(x_1,\ldots,x_m) = f_i(x_1,\ldots,x_m,a_1,\ldots,a_{l-m})\\
  \forall i \in \{1,\ldots,l-m\}&\,  \forall x_1,\ldots,x_m\in A:\nonumber\\
      &\, f_{n+i}(x_1,\ldots,x_m,a_1,\ldots,a_{l-m}) = a_i\\
  \forall b_1,\ldots,b_m\in A:&\, (x_1,\ldots,x_{l-m}) \mapsto (f_{n+1}(b_1,\ldots,b_m,x_1,\ldots,x_{l-m}),\ldots\\
  &\hspace{12mm}\ldots f_{k}(b_1,\ldots,b_m,x_1,\ldots,x_{l-m})) \in B_{l-m}(A)
 \end{align*}
\end{defn}

Let $F\subseteq M(A)$ be some mappings on $A$.
Then let $T(F)$ be the set of mappings realised with temporary storage by $F$, 
$T_S(F)$ the set of mappings realised with strong temporary storage by $F$.

These concepts and language are introduced based on Toffoli's  work. 
When we compose maps, this temporary storage
 allows us to use the garbage again in the subsequent function, as we 
know what it looks like. Thus the extra inputs and outputs are in some sense not garbage, as they are appropriately recycled and 
reduce waste information. In the engineering perspective, this reduces waste heat in implementations.

Let's look at an example.
Let $A=\Z_n$, $f(x,y) = (2x+y,xy)$. 
Then the  map $g(x)=2x$ is realised with temporary storage by $\{f\}$, since $g(x)=f_1(x,0)$ and $f_2(x,0)=0$ for all $x$.
However $g$ is not realised with strong temporary storage, because $y\mapsto f_2(0,y) \not\in B_1(A)$. 

\begin{lemma}
 Let $A$ be a finite set,  $F \subseteq B(A)$. Then $T_S(F) \subseteq T(F) \subseteq B(A)$.
\end{lemma}
Proof:
The first inclusion is clear, as the definition of strong temporary storage is stricter than the definition of temporary storage.

Let $g \in T(F)$. This means there is some $f \in C(F)$, $a_1,\dots,a_{l-m}\in A$ satisfying the requirements $(\ref{eqt1})$ and $(\ref{eqt2})$ above.
We know that $f\in B(A)$ because $F \subseteq B(A)$. Thus $f$ is balanced, so $l=k$ and thus $n=m$.
By $(\ref{eqt2})$, $f$ fixes $\{(x_1,\ldots,x_m,a_1\dots,a_{l-m}) \mid x_1,\dots,x_m\in A\}$ as a set.
Because $f$ is a bijection, $f$ is a permutation of $\{(x_1,\ldots,x_m,a_1\dots,a_{l-m}) \mid x_1,\dots,x_m\in A\}$.
By the definition of $g$ in $(\ref{eqt1})$, $g$ is a permutation of $A^m$,
so $g \in B(A)$.
Thus $T(F) \subseteq B(A)$ and we are done.
\hfill$\Box$

%
%
%
%
%

\section{The Functions of a Revclone}

In Toffoli's ``Fundamental Theorem,'' our Corollary \ref{corToffFund}, he shows that all maps and 
thus all functions can be realised by the full reversible \toffalg.
If we look at clones, this means that  the clone of all functions on $A$ can be realised, i.e.\ ${\mathcal O}^A \subseteq KS(B(A))$.

Let $A$ be a set. Let $F \subseteq M(A)$. We call $S_1(F)=\{\mu((1),f)\mid f\in F\} =\{f_1 \mid f \in F\}$ the \emph{function set} of $F$. 

\begin{thm}
\label{thmfunctions}
 Let $F$ be a \toffalg. Then $S_1(F)$ is a linear term algebra with projections.
\end{thm}
Proof: The operations $\tau$ and $\zeta$ act the same in the reversible \toffalg{}  and the linear term algebra.
Let $f,g\in F$, then
$\nabla (f_1) = (\bar \tau (i_1\oplus f))_1$.
Also $f_1 * g_1 = (f \circ_1 g)_1$. Thus $F$ being closed as a multiclone means that $S_1(F)$ is closed as a linear term algebra.

Moreover, the constants $\pi_\alpha \in F$ mean that all projections are also in $S_1(F)$,
so with $(1\,i) \in S_n$,
$\pi^n_i = (\pi_{(1\,i)})_1$.
\hfill$\Box$

This also applies for reversible \toffalgs. Note that we obtain the $ \nabla$ operation ``for free'' from the other operators.
In general we see the following.

\begin{cor}
 Let $F$ be a multiclone. Then $S_1(F)$ is a clone.
\end{cor}
Proof: A  multiclone is a \toffalg{}  closed under the $\Delta$ and $\nabla$ operations on multiclones.
By Theorem \ref{thmfunctions} we have $\nabla$ operating on $S_1(F)$.
For any $f_1 \in S_1(F)$, $\Delta (f_1) = (\Delta f)_1$ and we have clone closure.
\hfill$\Box$

This justifies the use of the expression ``clone'' in our expression multiclone.

Note that the mapping $F\mapsto \{f_1 \mid f \in F\}$ is not injective.
Let $A=\{0,1,2,3,4\}$,
 $C = \bigcup_n GL(n,\Z_5) \subseteq B(A)$ a reversible \toffalg, 
 \begin{align*}
  D=\bigcup_n \{-1,1\}SGL(n,\Z_5) = \{f\in C\mid \det f \in \{-1,1\}\}
 \end{align*}
 a subclone.
 $D$ is a proper sub reversible \toffalg{}  of $C$, by the determinant.
Both $GL(n,\Z_5)$ and $D\cap B_n(A)$ are transitive on $\Z_5^n \setminus \{(0,\dots,0)\}$ 
so both reversible \toffalgs{}  have function sets that are 
$\bigcup_n \{f:A^n\rightarrow A \mid f(x_1,\ldots,x_n)=\sum a_ix_i,\, a_i \in \Z_5,\, (a_1,\ldots,a_n)\neq 0\}$,
the clone of nonzero linear forms on $\Z_5$.

\begin{cor}
 Let $A$ be a finite set, let $F$ be a reversible \toffalg. 
 Then $g\in S_1(F)$, $g$ of arity $m$ implies that for all $ a \in A,\,\vert g^{-1}(a)\vert = \vert A\vert^{(m-1)}$.
\end{cor}
Proof: This is a corollary of Theorem \ref{thmToffFund}.
Let $g = \mu((1),f)$ for some $f\in F$, $f$ has arity $m$.
By requiring that the reversible mapping has the same arity, we require $r=m$ and $n=1$ in the statement of the theorem.
Because $r= \max(m,1+\lceil\log_{\vert A\vert} max_{a\in A^n} \vert g^ {-1} (a)\vert \rceil) = m$, we know that
$\log_{\vert A\vert} max_{a\in A^n} \vert g^ {-1} (a)\vert \leq m-1$.
However the average value of $\vert g^ {-1} (a)\vert$ is $\vert A \vert^{m-1}$ so the average value must be equal to the maximal value, 
so all  $\vert g^ {-1} (a)\vert$ are equal and we are done.
\hfill$\Box$

%

\section{Some results about Realisation}
\label{secgeneration}

Here we collect some results on reversible mappings. 
Much of this is based upon section 5 in \cite{toff80}, with the extension to Theorem \ref{ThmOdd}.



\begin{defn}
 Let $A$ be a set.
 Let $\alpha$ be a permutation of $A$, $o\in A$ some constant, $n\in \N$.
 For $n>1$, define $TG(n,\alpha,o) \in B_n(A)$ by 
 \begin{align*}
  TG(n,\alpha,o)(x_1,\ldots,x_n)_i &= x_i \hspace{10mm} i < n\\
  TG(n,\alpha,o)(x_1,\ldots,x_n)_n &= \begin{cases}
                                       \alpha(x_n) &\mbox{ if }  x_1=\ldots=x_{n-1}=o \\
                                       x_n &\mbox{ otherwise}
                                      \end{cases}
 \end{align*}
 Let $TG(1,\alpha,o)(x) = \alpha(x)$.
 Then $TG(n,\alpha,o) \in {B}_n$ is an invertible mapping, the \emph{Toffoli Gate} induced by $n$, $\alpha$ and $o$.
\end{defn}

In \cite{toff80},  this mapping is defined for $A=\Z_2$ with $\alpha$ swapping 0 and 1, $o=1$, written as $\theta^{(n)} = TG(n,(0\, 1),1)$.

\begin{defn}
Let $A$ be a set.
 A permutation $f \in {B}_n$ is \emph{elementary} if there exist distinct $x,y\in A^n$ such that $f(x)=y,f(y)=x$ 
 and all other elements of $A^n$ are fixed.
 An elementary permutation is \emph{atomic} if $x$ and $y$ only differ in one entry, i.e.\ $x_i=y_i$ for all $i$ except one.
\end{defn}

\begin{lemma}
\label{atomicrealisation}
 Let $A$ be a set.
 Let $n\in \N$, let $o\in A$. Let $f\in {B}_n(A)$ be atomic. 
 Then $f$ can be isomorpically realised by $\{TG(n,\alpha,o)\mid \alpha \in S_A\} \cup \{TG(1,\alpha,o)\mid \alpha \in S_A\}$.
\end{lemma}
Proof:
Suppose $f$ exchanges $x,y\in A^n$ which differ at position $i$.
Then $\pi_{(i,n)} \circ_n f \circ_n \pi_{(i,n)}$ is atomic, exchanging two vectors that differ at position $n$.
Wlog suppose $f$ is of this form.
Let $\alpha_i := (o,x_i) \in S_A$ be a transposition of $A$.
Let $\beta_i = TG(1,\alpha_i,o)$ for $i<n$, $\beta_n$ the identity.
Then $\beta = \beta_1\oplus\ldots\oplus\beta_n \in  {B}_n$ is an involution, generated by $TG(1,\alpha_i,o)$ as a reversible \toffalg.
Note that $\beta(x)_i = \beta(y)_i = o$ for $i <n$.

Note that $f(\beta(z)) = \beta(z)$ unless $\beta(z) = x$ respectively $\beta(z)=y$. 
But $\beta(z) = x$ (resp.\ $\beta(z)=y$) iff $z_i = o$ for all $i<n$ and $z_n=x_n$ (resp.\ $Z_n=y_n$).
So $\beta \bullet f \bullet \beta$ fixes all elements of $A^n$ except $(o,\ldots,o,x_n)$ and $(o,\ldots,o,y_n)$, which it exchanges.
Thus $\beta \bullet f \bullet \beta = TG(n,(x_n,y_n),o)$, 
so $f$ is in the reversible \toffalg{}  generated by $\{TG(n,\alpha,o)\mid \alpha \in S_A\} \cup \{TG(1,\alpha,o)\mid \alpha \in S_A\}$.
\hfill$\Box$

\begin{thm}[\cite{toff80} Thm 5.1]
\label{thmatomrel}
Let $A$ be a set.
 Any $f \in {B}_n(A)$ can be isomorphically realised by atomic permutations.
\end{thm}
Proof:
If $f$ is the identity, then $f=i_n$ and we are done. 
Suppose $f$ is nontrivial.
Any permutation can be written as a product of elementary permutations. So wlog, let $f$ be elementary,
exchanging $x,\,y\in A^n$.

Define a sequence $(a^{(i)}\in A^n \mid 1\leq i \leq n+1)$ such that  
 for all $k\in \{1,\ldots n+1\}$,  $(a^{(k)})_i = y_i$ if $i<k$, $(a^{(k)})_i = x_i$ otherwise.
Then $a^{(1)}=x$, $a^{(n+1)}=y$.
Note that $a^{(i)}=a^{(i+1)}$ or they differ in position $i$.
Thus $f_i \in {B}_n(A)$ exchanging $a^{(i)}$ and $a^{(i+1)}$ is an atomic permutation.
Now the permutation $f= f_1\bullet f_2 \bullet \ldots\bullet  f_{n-1} \bullet f_n \bullet f_{n-1} \bullet \ldots \bullet  f_1$ 
is a product of atomic permutations,
so we are done.
\hfill$\Box$

\begin{cor}
\label{corrBnGen}
 Let $A$ be a finite set, say $A=\{1,\dots,k\}$. Then ${B}_n(A)$ is realised by $\{TG(n,\alpha,1)\mid  \alpha \in \{(1,2),(1,\dots,k)\}\} \cup 
\{(1,2),(1,\dots,k)\}$.
\end{cor}
Proof: We remind ourselves that $TG(1,\alpha,1)=\alpha$ and thus see that the four gates listed here
generate $\{TG(n,\alpha,1)\mid \alpha \in S_A\} \cup \{TG(1,\alpha,1)\mid \alpha \in S_A\}$. 
By  Lemma \ref{atomicrealisation} this is enough to generate all atomic permutations and thus by  
Theorem \ref{thmatomrel},
enough to generate $B_n(A)$.
\hfill$\Box$

This is one possible correct form of the result hoped for in Conjecture 1 of \cite{yangetal}.
The conjecture can be demonstrated to be false using 
calculations in GAP \cite{GAP4} for $A$ of order 5 which can be understood as follows,
using a description kindly supplied by Erkko Lehtonen in a private communication.

\begin{conj}[Conjecture 1 of \cite{yangetal}]
Let $m\geq 3$, $n\geq 2$ be integers, $A=\Z_m$.
Then $B_n(A)$ is generated by the following permutations, where addition is modulo $m$.
\begin{itemize}
 \item Let $(ij) \in S_n$, $\sigma_{ij}=\pi_{(ij)}$
 \item Let $i\in \{1,\dots,n\}$, for all $x_1,\dots,x_n\in A$, let 
 $\nu_i(x_1,\dots,x_n) = (x_1,\dots,x_{i-1},x_i+1,x_{i+1},\dots,x_n)$
 \item For all $x_1,\dots,x_n\in A$, let
 \[ \Upsilon(x_1,\dots,x_n) = \begin{cases}
                           (x_1+1,0,\dots,0) & \mbox{ if }x_2=\dots=x_n=0 \\
                           (x_1,\dots,x_n) & \mbox{otherwise}
                          \end{cases}
 \]
\end{itemize}
\end{conj}
This can be falsified by considering the parities of the permutations.

The permutation $\sigma_{ij}$ is a product of $m^{n-1}(m-1)/2$ transpositions.
Thus $\sigma_{ij}$ is odd iff $m \equiv 3 \mod 4$ or ($m \equiv 2 \mod 4$ and $n=2$).

Every element of $A^n$ belongs to a cycle of length $m$ of $\nu_i$.
Thus $\nu_i$ is a product of $m^{n-1}$ cycles of length $m$, each of which is
made of $m-1$ transpositions, so $\nu_i$ is a product of $m^{n-1}(m-1)$ transpositions,
which is always even.

The permutation $\Upsilon$ is a cycle of length $m$, thus is even iff $m$ is odd.

Thus for $m=5$ or more generally $m \equiv 1 \mod 4$ all of these permutations are even so they cannot generate $B_n(A)$ for any $n$.

For $m$ even, $n\geq 2$, the conjecture might hold,  calculations in GAP for small values find no contradiction.
So we obtain the following conjecture that a smaller generating set might suffice.
\begin{conj}
 For $A=\{1,\dots,m\}$ even, $n\geq 2$, \[B_n(A) = C(\{(1,\dots,m), TG(n,(1,\dots,m),1)\})\]
\end{conj}

Below we will see that a small generating set exists for $A$ of odd order.

It is worth noting that the main result in \cite{yangetal} is also subtly wrong for $n=1$, as  the permutation $(0\,1)$ is not generated.
The authors use an interesting property of ternary logic 
that any permutation of $\{0,1,2\}$ (written in $\Z_3$) can be written as $\pmtv[0 1 2] {{0}{n} {1}{n+1} {2}{n+1+1} }$
or $\pmtv[0 1 2] {{0}{n} {1}{n+2} {2}{n+2+2} }$ for some $n$, 
leading to their minimal generating set for $n\geq 2$.
For this reason it might be that their result only applies for ternary alphabets.

The following results lets us calculate precisely what the $n$-ary part of a reversible \toffalg{}  looks like, given
the generators.

\begin {thm}
\label{thmFbar}
  Let $A$ be a set.
  Let $F \subseteq B(A)$, $n\in \N$. 
  Let 
  \begin{align*}
     \bar F = &\{ \bar f \mid f \in F,\,\mbox{ arity } f =m < n,\, \bar f = f\oplus i_{n-m}\}\\
       &\cup\{f \in F\mid \mbox{ arity }f=n\}
         \cup \{\pi_\alpha \mid \alpha \in S_n\}.
  \end{align*}

Then the group $(\langle \bar F \rangle ; \circ_n)$ is equal to $(C(F) \cap B_n(A); \circ_n)$.
\end {thm}
Proof:
Both groups are subgroups of $(B_n(A);\circ_n)$, so we need only prove they are equal as sets.

We proceed by induction.
For the case $n=1$, $\bar F$ consists of the permutations of $A$ within $F$.
Then the subgroup of $S_A$ generated by these permutations is precisely $C(F) \cap B_1(A) = C(F) \cap S_A$.
Then we proceed with our induction.
Suppose our claim holds up to $n-1$.

$(\subseteq):$ Every element $\bar f = f \oplus i_{n-m}$ of $\bar F$ is within $C(F)$. 
The arity of all elements of $\bar F$ is $n$ so they are all within $B_n(A)$. So this inclusion holds.

$(\supseteq):$ 
Let $f\in C(F) \cap B_n(A)$, so $f$ is either in $\bar F$ or is a term of the form:
\begin{enumerate}
 \item $f = g \oplus h$  for some $g,h \in  C(F)$, both of smaller arity than $n$.
 \item $f = g \circ_k h$ for some $g,h \in  C(F)$, both of smaller arity than $n$, $k<n$.
 \item $f = g \circ_k h$ for some $g,h \in  C(F)$, exactly one of them in $B_n(A)$, so the other one has arity $k$, $k<n$.
 \item $f = g \circ_n h$ for some $g,h \in  C(F)\cap B_n(A)$.
\end{enumerate}
In case 1, let $m$ be the arity of $g$, so $h$ has arity $n-m$. 
We know from our induction hypothesis that $g$ is a product $g= \bar g_1\circ_m\ldots\circ_m\bar g_l$ of 
elements of the form $\bar g_j = f_j\oplus i_{m_j}$ for some $f_j\in F$ with arity less than $m$,  
$\bar g_j = f_j$ for some $f_j\in F$ with arity  $m$,
or $\bar g_j = \pi_{\alpha_j}$ for some 
permutation $\alpha_j\in S_m$.
Then let $\phi_j = f_j \oplus i_{m_j + n-m} =  f_j \oplus i_{m_j} \oplus i_{n-m}\in \bar F$ respectively
$\phi_j = f_j \oplus i_{ n-m}\in \bar F$ respectively
$\phi_j = \pi_{\alpha_j}\oplus i_{n-m} = \pi_{\beta_j}$ where $\beta_j \in S_n$
is equal to $\alpha_j$ on $\{1,\ldots,m\}$ and
fixes the elements $\{m+1,\ldots,n\}$.
Then $\phi_1\circ_n\ldots \circ_n \phi_l = g\oplus i_{n-m}$, so $g\oplus i_{n-m} \in \langle \bar F \rangle$.

Similarly we can write $h\oplus i_m$ as an element of $\langle \bar F \rangle$.
Let $\delta$ be the permutation of $\{1,\ldots,n\}$ defined by
\[\delta = \pmtv[1\dots {m}{(m+1)}\dots {n}] {{1}{(n-m+1)} {m}{n} {(m+1)}{1} {n}{(n-m)} }.\]
Then 
\[\pi_{\delta^{-1}} \circ_n (h\oplus i_m) \circ_n \pi_\delta \circ_n (g\oplus i_{n-m}) = g \oplus h\]
Thus we see that $f \in (\langle \bar F \rangle ; \circ_n)$.

In case 2 we proceed similarly. 
In this case $g$ is $m$-ary and $h$ is $n-m+k$-ary.
We use the same techniques as above to show that $g\oplus i_{n-m}, h\oplus i_{m-k} \in \langle \bar F \rangle$.
Let $\delta$ be the permutation of $\{1,\ldots,n\}$ given by
\[\delta = \pmtv[1\ldots {k}{(k+1)}\ldots {(n-m+k)}{(n-m+k+1)}\ldots {n}] {{1}{1} {k}{k} {(k+1)}{(m+1)} {(n-m+k)}{n}{(n-m+k+1)}{(k+1)} {n}{m} }\]
so $\pi_\delta \in \bar F$.
Then 
\begin{align*}
 &(g\oplus i_{n-m}) \circ_n \pi_\delta \circ_n (h\oplus i_{m-k})(x_1,\ldots,x_n)\\
 & = (g\oplus i_{n-m}) \circ_n \pi_\delta (h_1(x_1,\ldots,x_{n-m+k}),\ldots,h_{n-m+k}(x_1,\ldots,x_{n-m+k}),\\
      & \hspace{30mm}x_{n-m+k+1},\dots,x_n)\\
 & = (g\oplus i_{n-m})(h_1(x_1,\ldots,x_{n-m+k}),\ldots,h_{k}(x_1,\ldots,x_{n-m+k}),\\
    & \hspace{30mm}x_{n-m+k+1},\dots,x_n,\\
                            & \hspace{30mm} h_{k+1}(x_1,\ldots,x_{n-m+k}),\dots,h_{n-m+k}(x_1,\ldots,x_{n-m+k})) \\
 &= (g_1(h_1(x_1,\ldots,x_{n-m+k}),\ldots,h_{k}(x_1,\ldots,x_{n-m+k}),x_{n-m+k+1},\dots,x_n),\\
 & \hspace {20mm} \vdots\\
 & \hspace {6mm} g_m(h_1(x_1,\ldots,x_{n-m+k}),\ldots,h_{k}(x_1,\ldots,x_{n-m+k}),x_{n-m+k+1},\dots,x_n),\\
 & \hspace {6mm} h_{k+1}(x_1,\ldots,x_{n-m+k}),\dots,h_{n-m+k}(x_1,\ldots,x_{n-m+k})) \\
 & = g \circ_k h (x_1,\dots,x_n)
\end{align*}
Thus we see that $f \in (\langle \bar F \rangle ; \circ_n)$.

In case 3 we have  $k < n$ then one of the terms is of lower arity, 
so as above we can write $g \oplus i_{n-k}$ (respectively $h\oplus i_{n-k}$)
as an element of $B_n(A)$ and have $f = (g \oplus i_{n-k}) \circ_n h$ (respectively $f=g \circ_n (h\oplus i_{n-k})$).
Then we have the next case.

In case 4  we have two terms of arity $n$ but of strictly lower term complexity.
We use induction on term complexity. 
For the initial case, we look at the trivial terms  $t \in C(F)\cap B_n(A)$.
These are either elements of $F$ or the permutations $\pi_\alpha$ for $\alpha \in S_n$.
In both these cases, $t \in \bar F$ so $t\in \langle \bar F \rangle$.
Now we look at $f = g \circ_n h$. Both $g$ and $h$ have lower term complexity than $f$.
Thus we know that $g,h \in \langle \bar F \rangle$.
Thus $f  = g \circ_n h\in  \langle \bar F \rangle$ and we are done.
\hfill$\Box$

We can apply this to obtain a generalisation of one of Toffoli's results.

\begin{cor}[\cite{toff80} Thm 5.2]
\label{thmtoffeven}
 Let $A$ be a finite set of even order, $n\in \N$.
 Then ${B}_n(A)$ is not isomorphically realised by $\{TG(i,\alpha,o) \mid \alpha \in S_A,\, i<n\}$ for any $o\in A$
\end{cor}
Proof:
For $n=1$ the set of generators is empty, so we only have the identity permutation.
For $n=2$, we calculate the value of $\vert C(\{TG(i,\alpha,o) \mid \alpha \in S_A,\, i<n\}) \cap B_2(A)\vert$ and
$\vert B_2(A)\vert$.
First note that $\{TG(i,\alpha,o) \mid \alpha \in S_A,\, i<n\} = \{TG(1,\alpha,o) \mid \alpha \in S_A\}$
so our generating set consists only of permutations of entries in $A$.
We see that 
$C(\{TG(i,\alpha,o) \mid \alpha \in S_A,\, i<n\}) \cap B_2(A)$
thus consists only of bijections that map $(a,b)\in A^2$ to $(a^\alpha,b^\beta)$ or $(b^\beta, a^\alpha)$ for
some $\alpha,\beta \in S_A$.
Thus  
$\vert C(\{TG(i,\alpha,o) \mid \alpha \in S_A,\, i<2\}) \cap B_2(A)\vert = \vert S_A\vert^2 2=2 (\vert A\vert !)^2$.
On the other hand, $\vert B_2(A)\vert = (\vert A \vert^2)!$ which is strictly greater than $2 (\vert A\vert !)^2$ for all $A$
with at least two elements.

We carry on for the other cases.
Using the terminology in Theorem \ref{thmFbar}, we note that $\overline{TG}(i,\alpha,o) = TG(i,\alpha,o) \oplus i_{n-i}$.
The action of $\tau = TG(i,\alpha,o) \oplus i_{n-i}$ on $A^n$, when not identity, is of the form
\[\tau (o,\ldots,o,a,a_1,\ldots,a_{n-i}) = (o,\ldots,o,\alpha(a),a_1,\ldots,a_{n-i}) \]
If we write $\alpha$ as a product of transpositions $\alpha = \alpha_1\ldots\alpha_k$ 
we can use the expression above to see that there will be
$k \vert A\vert^{n-i}$ transpositions when we write the action of $\tau$ on $A^n$.
Thus $\overline{TG}(i,\alpha,o) = TG(i,\alpha,o) \oplus i_{n-i}$ is an even permutation in $B_n(A)$.

The action of $\pi_\beta$ for $\beta \in S_n$ acting on $A^n$ can also be calculated. 
If $\beta$ is a transposition in $S_n$, say $\beta = (ij)$, 
then $\pi_\beta$   acts nontrivially on $(a_1,\ldots, a_i, \ldots, a_j, \ldots, a_n)$ with $a_i \neq a_j$. There are 
$\frac{\vert A\vert (\vert A \vert -1)}{2} (\vert A\vert )^{n-2}$ such tuple pairs. This number is
even when $\vert A\vert$ is even, except in the case $n= 2$ (when $\vert A \vert\equiv 2 \mod 4$) 
which we dealt with above. 

From Theorem \ref{thmFbar}, we know that the arity $n$ part of $C(\{TG(i,\alpha,o) \mid \alpha \in S_A,\, i<n\})$ 
is generated
by precisely these permutations, which are all even. 
Thus  for some $a,b\in A$, $a\neq b$, $TG(n,(a,b),a)$ is not in the generated reversible \toffalg{}
and thus $B_n(A)$ is not isomorphically realised by $\{TG(i,\alpha,o) \mid \alpha \in S_A,\, i<n\}$.
\hfill$\Box$

The following has been noted using examples in GAP\cite{GAP4}.
\begin{conj}
 Let $A$ be a finite set of even order, $n\geq 3$,  $F=\{TG(i,\alpha,o) \mid \alpha \in S_A,\, i<n, o\in A\}$.
 We claim that $C(F) \cap B_n$ is a subgroup of $B_n$ isomorphic to
 the alternating group on $A^n$, except for $\vert A \vert = 2$, $n=3$.
\end{conj}

{ \it Note added in referee process of this manuscript: The conjecture has been proven in
\cite{boykettkariville17}.}

The restriction shown in Corollary \ref{thmtoffeven}  does not hold for $A$ odd, as we see in the following result. 
Moreover, it shows that we only need use the Toffoli gates of arity 1 and 2 to obtain all bijections on $A$.

\begin{thm}
\label{ThmOdd}
 Let $A $ be a finite set of odd order.
 Then $B(A)$  is isomorphically realised by 
 \[T_2=\{TG(i,\alpha,o) \mid \alpha \in S_A,\, i<3\} = S_A \cup \{TG(2,\alpha,o) \mid \alpha \in S_A\}\] 
 for any $o\in A$.
\end{thm}
Proof:
Let $\vert A \vert = k$, wlog $A=\{1,\ldots,k\}$ and $o=1$. 
We proceed by induction on $n$, i.e.\ showing that $B_n(A) \subseteq C(T_2)$. Our start is for $n=2$, given by the hypothesis.
We assume we have shown the claim up to $n$.
We first show that we can obtain $TG(n+1,(1 2),o)$ from $\{TG(i,\alpha,o) \mid \alpha \in S_A,\, i\leq n\}$.
Let $\gamma =(n\; n+1) \in S_{n+1}$.


Define
\begin{align*}
 \Sigma_1 = (TG&(n,(1\ldots k)^{-1},1) \oplus i_1)  \bullet \pi_{\gamma}  \bullet (TG(n,(1 2),1) \oplus i_1)\\  
   &\bullet \pi_{\gamma} \bullet (i_{n-1} \oplus TG(2,(1 2),1) )\\
   &\bullet (TG(n,(1\ldots k),1) \oplus i_1) \in B_{n+1}(A)
\end{align*}

We calculate.

\begin{align*}
   &\Sigma_1(x_1,\ldots,x_{n+1}) \\
 = & \begin{cases}
      (TG(n,(1\ldots k)^{-1},1) \oplus i_1)  \bullet \ldots& \bullet (i_{n-1} \oplus TG(2,(1 2),1) )(x_1,\ldots,x_{n-1},x_{n}^{(1\ldots k)},x_{n+1})\\ & \mbox{ if } x_1=\ldots=x_{n-1}=1 \\
      (TG(n,(1\ldots k)^{-1},1) \oplus i_1)  \bullet \dots& \bullet (i_{n-1} \oplus TG(2,(1 2),1) )(x_1,\ldots,x_{n-1},x_{n},x_{n+1})\\ & \mbox{ otherwise}
     \end{cases}
\end{align*}
\begin{align*} 
 = & \begin{cases}
      (TG(n,(1\ldots k)^{-1},1) \oplus i_1)  \bullet \dots& \bullet \pi_{\gamma}(x_1,\ldots,x_{n-1},x_{n}^{(1\ldots k)},x_{n+1}^{(12)})\\ & \hspace{8mm}\mbox{ if } x_1=\ldots=x_{n-1}=1, x_n=k \\
      (TG(n,(1\ldots k)^{-1},1) \oplus i_1)  \bullet \dots& \bullet \pi_{\gamma}(x_1,\ldots,x_{n-1},x_{n}^{(1\ldots k)},x_{n+1})\\ & \hspace{8mm}\mbox{ if } x_1=\ldots=x_{n-1}=1, x_n\neq k \\
      (TG(n,(1\ldots k)^{-1},1) \oplus i_1)  \bullet \dots& \bullet \pi_{\gamma}(x_1,\ldots,x_{n-1},x_{n},x_{n+1}^{(12)})\\ & \hspace{8mm}\mbox{ if some } x_1,\ldots,x_{n-1} \neq 1, x_n=1\\
      (TG(n,(1\ldots k)^{-1},1) \oplus i_1)  \bullet \dots& \bullet \pi_{\gamma}(x_1,\ldots,x_{n-1},x_{n},x_{n+1})\\ & \hspace{8mm}\mbox{ if some } x_1,\ldots,x_{n-1} \neq 1, x_n\neq 1
     \end{cases}   \\
 = &\begin{cases}
      (TG(n,(1\ldots k)^{-1},1) \oplus i_1)  \bullet  \pi_{\gamma} &\bullet (TG(n,(1 2),1) \oplus i_1)(x_1,\ldots,x_{n-1},x_{n+1}^{(12)},x_{n}^{(1\ldots k)})\\ & \hspace{8mm}\mbox{ if } x_1=\ldots=x_{n-1}=1, x_n=k \\
      (TG(n,(1\ldots k)^{-1},1) \oplus i_1)  \bullet  \pi_{\gamma} &\bullet (TG(n,(1 2),1) \oplus i_1)(x_1,\ldots,x_{n-1},x_{n+1},x_{n}^{(1\ldots k)})\\ & \hspace{8mm}\mbox{ if } x_1=\ldots=x_{n-1}=1, x_n\neq k \\
      (TG(n,(1\ldots k)^{-1},1) \oplus i_1)  \bullet  \pi_{\gamma} &\bullet (TG(n,(1 2),1) \oplus i_1)(x_1,\ldots,x_{n-1},x_{n+1}^{(12)},x_{n})\\ & \hspace{8mm}\mbox{ if some } x_1,\ldots,x_{n-1} \neq 1, x_n=1\\
      (TG(n,(1\ldots k)^{-1},1) \oplus i_1)  \bullet  \pi_{\gamma} &\bullet (TG(n,(1 2),1) \oplus i_1)(x_1,\ldots,x_{n-1},x_{n+1},x_{n})\\ &\hspace{8mm} \mbox{ if some } x_1,\ldots,x_{n-1} \neq 1, x_n\neq 1
     \end{cases}\\
 = &\begin{cases}
      (TG(n,(1\ldots k)^{-1},1) \oplus i_1)  \bullet &\pi_{\gamma}(x_1,\ldots,x_{n-1},x_{n+1}^{(12)(12)},x_{n}^{(1\ldots k)})\\ & \hspace{8mm}\mbox{ if } x_1=\ldots=x_{n-1}=1,x_n=k \\
      (TG(n,(1\ldots k)^{-1},1) \oplus i_1)  \bullet &\pi_{\gamma}(x_1,\ldots,x_{n-1},x_{n+1}^{(12)},x_{n}^{(1\ldots k)})\\ & \hspace{8mm}\mbox{ if } x_1=\ldots=x_{n-1}=1, x_n\neq k \\
      (TG(n,(1\ldots k)^{-1},1) \oplus i_1)  \bullet &\pi_{\gamma}(x_1,\ldots,x_{n-1},x_{n+1}^{(12)},x_{n})\\ & \hspace{8mm}\mbox{ if some } x_1,\ldots,x_{n-1} \neq 1, x_n=1\\
      (TG(n,(1\ldots k)^{-1},1) \oplus i_1)  \bullet &\pi_{\gamma}(x_1,\ldots,x_{n-1},x_{n+1},x_{n})\\ & \hspace{8mm}\mbox{ if some } x_1,\ldots,x_{n-1} \neq 1, x_n\neq 1
     \end{cases}
\end{align*}

\begin{align*}
 = &\begin{cases}
      (TG(n,(1\ldots k)^{-1},1) \oplus i_1) &(x_1,\ldots,x_{n-1},x_{n}^{(1\ldots k)},x_{n+1}) \\ & \hspace{8mm} \mbox{ if } x_1=\ldots=x_{n-1}=1, x_n=k \\
      (TG(n,(1\ldots k)^{-1},1) \oplus i_1) &(x_1,\ldots,x_{n-1},x_{n}^{(1\ldots k)},x_{n+1}^{(12)}) \\ & \hspace{8mm} \mbox{ if } x_1=\ldots=x_{n-1}=1, x_n\neq k \\
      (TG(n,(1\ldots k)^{-1},1) \oplus i_1) &(x_1,\ldots,x_{n-1},x_{n},x_{n+1}^{(12)}) \\ & \hspace{8mm} \mbox{ if some } x_1,\ldots,x_{n-1} \neq 1, x_n=1\\
      (TG(n,(1\ldots k)^{-1},1) \oplus i_1) &(x_1,\ldots,x_{n-1},x_{n},x_{n+1}) \\ & \hspace{8mm} \mbox{ if some } x_1,\ldots,x_{n-1} \neq 1, x_n\neq 1
     \end{cases}\\
 = &\begin{cases}
      (x_1,\ldots,x_{n-1},x_{n},x_{n+1}) & \mbox{ if } x_1=\ldots=x_{n-1}=1, x_n=k \\
      (x_1,\ldots,x_{n-1},x_{n},x_{n+1}^{(12)}) & \mbox{ if } x_1=\ldots=x_{n-1}=1, x_n\neq k \\
      (x_1,\ldots,x_{n-1},x_{n},x_{n+1}^{(12)}) & \mbox{ if some } x_1,\ldots,x_{n-1} \neq 1, x_n=1\\
      (x_1,\ldots,x_{n-1},x_{n},x_{n+1}) & \mbox{ if some } x_1,\ldots,x_{n-1} \neq 1, x_n\neq 1
     \end{cases}
\end{align*}

We see that $\Sigma_1$ is nonidentity iff:
\begin{enumerate}
 \item Not all of $x_1\ldots,x_{n-1}$ are 1 and $x_n=1$, then $x_{n+1} \mapsto x_{n+1}^{(12)}$, or
 \item $x_1=\ldots=x_{n-1}=1$ and $x_n\neq k$, then $x_{n+1} \mapsto x_{n+1}^{(12)}$.
\end{enumerate}

For $m \in \{2,\ldots,k-1\}$, let
\begin{align*}
 \sigma_m = (TG(n,(1 m),1) \oplus i_1)  \bullet (i_{n-1} \oplus TG(2,(1 2),1)) \bullet(TG(n,(1 m),1) \oplus i_1) 
 \in B_{n+1}(A)
\end{align*}

We calculate 
\begin{align*}
 &\sigma_m(x_1,\ldots,x_{n+1}) \\
 &=  \begin {cases}
      (TG(n,(1 m),1) \oplus i_1)  \bullet (i_{n-1} \oplus TG(2,(1 2),1))(x_1,\ldots,x_{n-1},x_n^{(1m)},x_{n+1}) \\ \hspace{58mm}\mbox{ if } x_1=\ldots=x_{n-1}=1 \\
      (TG(n,(1 m),1) \oplus i_1)  \bullet (i_{n-1} \oplus TG(2,(1 2),1))(x_1,\ldots,x_{n-1},x_n,x_{n+1}) \\ \hspace{58mm}\mbox{ if some } x_1,\ldots, x_{n-1}\neq 1
     \end {cases} 
 \end{align*}    
 \begin{align*}    
 &=  \begin {cases}
      (TG(n,(1 m),1) \oplus i_1)  (x_1,\ldots,&x_n^{(1m)},x_{n+1}^{(12)}) \\&\mbox{ if } x_1=\ldots=x_{n-1}=1, x_n=m\\
      (TG(n,(1 m),1) \oplus i_1)  (x_1,\ldots,&x_n^{(1m)},x_{n+1}) \\&\mbox{ if } x_1=\ldots=x_{n-1}=1, x_n\neq m\\
      (TG(n,(1 m),1) \oplus i_1)  (x_1,\ldots,&x_n,x_{n+1}^{(12)}) \\&\mbox{ if some } x_1,\ldots, x_{n-1}\neq 1, x_n=1\\
      (TG(n,(1 m),1) \oplus i_1)  (x_1,\ldots,&x_n,x_{n+1}) \\&\mbox{ if some } x_1,\ldots, x_{n-1}\neq 1, x_{n} \neq 1
     \end {cases} \\
& =  \begin {cases}
        (x_1,\ldots,x_n,x_{n+1}^{(12)}) &\mbox{ if } x_1=\ldots=x_{n-1}=1, x_n=m\\
        (x_1,\ldots,x_n,x_{n+1}) &\mbox{ if } x_1=\ldots=x_{n-1}=1, x_n\neq m\\
        (x_1,\ldots,x_n,x_{n+1}^{(12)}) &\mbox{ if some } x_1,\ldots, x_{n-1}\neq 1, x_n=1\\
        (x_1,\ldots,x_n,x_{n+1}) &\mbox{ if some } x_1,\ldots, x_{n-1}\neq 1, x_{n} \neq 1
     \end {cases}
\end{align*}

Once again this function is almost always identity, $\sigma_m$ is nonidentity iff:
\begin{enumerate}
 \item Not all of $x_1\ldots,x_{n-1}$ are 1 and $x_n=1$, then $x_{n+1} \mapsto x_{n+1}^{(12)}$, or
 \item $x_1=\ldots=x_{n-1}=1$ and $x_n= m$, then $x_{n+1} \mapsto x_{n+1}^{(12)}$.
\end{enumerate}

Then define
\begin{align*}
 \Sigma_2 = \sigma_{k-1} \bullet \dots \bullet\sigma_{2}
\end{align*}
We have then that $\Sigma_2$ is nonidentity iff
\begin{enumerate}
 \item Not all of $x_1,\ldots,x_{n-1}$ are 1 and $x_n=1$, then $x_{n+1} \mapsto x_{n+1}^{(12)^{k-2}}=x_{n+1}^{(12)}$ because $k$ and thus $k-2$ are odd, or
 \item $x_1=\ldots=x_{n-1}=1$ and $x_n \in \{2,\ldots,k-1\}$, then $x_{n+1} \mapsto x_{n+1}^{(12)}$.
\end{enumerate}

We see that $\Sigma_2 \bullet \Sigma_1$ is identity unless one of the factors is nonidentity. There are three cases:
\begin{enumerate}
 \item Both are nonidentity by their first case. Then not all of $x_1,\ldots,x_{n-1}$ are 1 and $x_n=1$, so
 \begin{align*}
  \Sigma_2 \bullet \Sigma_1 (x_1\ldots x_{n+1}) 
  & = \Sigma_2(x_1,\ldots,x_{n+1}^{(12)}) \\
  & = (x_1,\ldots,x_{n+1}^{(12)(12)}) \\
  & = (x_1,\ldots,x_{n+1})
 \end{align*}
 so $\Sigma_2 \bullet \Sigma_1$ is the identity.
 \item Both are nonidentity by their second case, so $x_1=\ldots=x_{n-1}=1$ and $x_n \in \{2,\ldots,k-1\}$, 
 \begin{align*}
  \Sigma_2 \bullet \Sigma_1 (x_1\ldots x_{n+1}) 
  & = \Sigma_2(x_1,\ldots,x_{n+1}^{(12)}) \\
  & = (x_1,\ldots,x_{n+1}^{(12)(12)}) \\
  & = (x_1,\ldots,x_{n+1})
 \end{align*}
 so $\Sigma_2 \bullet \Sigma_1$ is the identity.
 \item Only $\Sigma_1$ is nonidentity by the second case, so $x_1=\ldots=x_{n}=1$, then 
 \begin{align*}
  \Sigma_2 \bullet \Sigma_1 (x_1\ldots x_{n+1}) 
  & = \Sigma_2(x_1,\ldots,x_{n+1}^{(12)}) \\
  & = (x_1,\ldots,x_{n+1}^{(12)}) \\
  & = (x_1,\ldots,x_{n+1}^{(12)})
 \end{align*}
 so we have the only case that $\Sigma_2 \bullet\Sigma_1$ is nonidentity.
\end{enumerate}
Thus we see that $\Sigma_2 \bullet\Sigma_1 = TG(n+1,(1 2),1)$.

We now look at $TG(n+1,(1 \ldots k),1)$.
Because the group generated by $(1 \ldots k)$ is cyclic of odd order, the homomorphism $x\mapsto x^2$ of this group is an automorphism.
Thus there exists some $\beta$ such that $\beta^2 = (1 \ldots k)$. This will also be a $k$-cycle, write $\beta = (\beta_1 \ldots \beta_k)$.
Let $k=2l+1$. 
Define $\alpha = (\beta_1 \beta_k)(\beta_2 \beta_{k-1}) \ldots (\beta_l \beta_{l+2})$.
Then $\alpha \beta \alpha = \beta^{-1}$, so $\beta\alpha\beta^{-1}\alpha=\beta^2=(1\ldots k)$.

Let $\gamma = (n\,n+1) \in S_{n+1}$. Define 
\begin{align*}
\begin{split}
 \Sigma = \pi_\gamma \bullet (TG(n,\alpha,1) \oplus i_1)\bullet\pi_\gamma \bullet (i_{n-1} \oplus TG(2,\beta^{-1},1)) 
         \bullet\pi_\gamma \hspace{10mm}\\
 \bullet(TG(n,\alpha,1) \oplus i_1)\bullet\pi_\gamma\bullet (i_{n-1} \oplus TG(2,\beta,1)) \in B_{n+1}(A)
 \end{split}
\end{align*}
We calculate.
\begin{align*}
 &\Sigma(x_1,\ldots,x_{n+1}) \\
 &= \begin{cases}
     \pi_\gamma \bullet \dots \bullet \pi_\gamma (x_1,\ldots,x_n,x_{n+1}^\beta) &\mbox{ if } x_n=1 \\
     \pi_\gamma \bullet \dots \bullet \pi_\gamma (x_1,\ldots,x_n,x_{n+1}) & \mbox{ if }x_n\neq 1 \\
    \end{cases} \\
 &= \begin{cases}
     \pi_\gamma \bullet \dots \bullet (TG(n,\alpha,1) \oplus i_1) (x_1,\ldots,x_{n+1}^\beta,x_n) &\mbox{ if } x_n=1 \\
     \pi_\gamma \bullet \dots \bullet (TG(n,\alpha,1) \oplus i_1) (x_1,\ldots,x_{n+1},x_n) & \mbox{ if }x_n\neq 1 
    \end{cases} \\
 &= \begin{cases}
     \pi_\gamma \bullet \dots \bullet \pi_\gamma (x_1,\ldots,x_{n+1}^{\beta\alpha},x_n) &\mbox{ if } x_1=\ldots=x_n=1 \\
     \pi_\gamma \bullet \dots  \bullet\pi_\gamma (x_1,\ldots,x_{n+1}^\beta,x_n) & \mbox{ if some } x_1,\ldots,x_{n-1} \neq 1,\,x_n=1 \\
     \pi_\gamma \bullet \dots \bullet \pi_\gamma (x_1,\ldots,x_{n+1}^\alpha,x_n) & \mbox{ if }x_1=\ldots=x_{n-1}=1,\, x_n\neq 1 \\
     \pi_\gamma \bullet \dots \bullet \pi_\gamma (x_1,\ldots,x_{n+1},x_n) & \mbox{ if some } x_1,\ldots,x_{n-1} \neq 1,\,x_n\neq 1 
    \end{cases}  
\end{align*} 
\begin{align*}   
 &= \begin{cases}
     \pi_\gamma \bullet  \dots \bullet (i_{n-1} \oplus TG(2,\beta^{-1},1)) &(x_1,\ldots,x_n,x_{n+1}^{\beta\alpha}) \\& \mbox{ if }x_1=\ldots=x_n=1 \\
     \pi_\gamma \bullet  \dots \bullet (i_{n-1} \oplus TG(2,\beta^{-1},1)) &(x_1,\ldots,x_n,x_{n+1}^\beta) \\& \mbox{ if some } x_1,\ldots,x_{n-1} \neq 1,\,x_n=1 \\
     \pi_\gamma \bullet  \dots \bullet (i_{n-1} \oplus TG(2,\beta^{-1},1)) &(x_1,\ldots,x_n,x_{n+1}^\alpha) \\&\mbox{ if } x_1=\ldots=x_{n-1}=1,\, x_n\neq 1 \\
     \pi_\gamma \bullet  \dots \bullet (i_{n-1} \oplus TG(2,\beta^{-1},1)) &(x_1,\ldots,x_n,x_{n+1}) \\& \mbox{ if some } x_1,\ldots,x_{n-1} \neq 1,\,x_n\neq 1 
    \end{cases}\\
 &= \begin{cases}
     \pi_\gamma \bullet  \dots \bullet \pi_\gamma (x_1,\ldots,x_n,x_{n+1}^{\beta\alpha\beta^{-1}}) & \mbox{ if }x_1=\ldots=x_n=1 \\
     \pi_\gamma \bullet  \dots \bullet \pi_\gamma (x_1,\ldots,x_n,x_{n+1}^{\beta\beta^{-1}}) & \mbox{ if some } x_1,\ldots,x_{n-1} \neq 1,\,x_n=1 \\
     \pi_\gamma \bullet  \dots \bullet \pi_\gamma (x_1,\ldots,x_n,x_{n+1}^\alpha) &\mbox{ if } x_1=\ldots=x_{n-1}=1,\, x_n\neq 1 \\
     \pi_\gamma \bullet  \dots \bullet \pi_\gamma (x_1,\ldots,x_n,x_{n+1}) & \mbox{ if some } x_1,\ldots,x_{n-1} \neq 1,\,x_n\neq 1 
    \end{cases}\\
 &= \begin{cases}
     \pi_\gamma \bullet (TG(n,\alpha,1) \oplus i_1) (x_1,\ldots,x_{n+1}^{\beta\alpha\beta^{-1}},x_n) &\mbox{ if } x_1=\ldots=x_n=1 \\
     \pi_\gamma \bullet (TG(n,\alpha,1) \oplus i_1) (x_1,\ldots,x_{n+1},x_n) & \mbox{ if some } x_1,\ldots,x_{n-1} \neq 1,\,x_n=1 \\
     \pi_\gamma \bullet (TG(n,\alpha,1) \oplus i_1) (x_1,\ldots,x_{n+1}^\alpha,x_n) &\mbox{ if } x_1=\ldots=x_{n-1}=1,\, x_n\neq 1 \\
     \pi_\gamma \bullet (TG(n,\alpha,1) \oplus i_1) (x_1,\ldots,x_{n+1},x_n) & \mbox{ if some } x_1,\ldots,x_{n-1} \neq 1,\,x_n\neq 1 
    \end{cases}\\
 &= \begin{cases}
     \pi_\gamma  (x_1,\ldots,x_{n+1}^{\beta\alpha\beta^{-1}\alpha},x_n) & \mbox{ if }x_1=\ldots=x_n=1 \\
     \pi_\gamma  (x_1,\ldots,x_{n+1},x_n) & \mbox{ if some } x_1,\ldots,x_{n-1} \neq 1,\,x_n=1 \\
     \pi_\gamma  (x_1,\ldots,x_{n+1}^{\alpha\alpha},x_n) &\mbox{ if } x_1=\ldots=x_{n-1}=1,\, x_n\neq 1 \\
     \pi_\gamma  (x_1,\ldots,x_{n+1},x_n) & \mbox{ if some } x_1,\ldots,x_{n-1} \neq 1,\,x_n\neq 1 
    \end{cases}\\
 &= \begin{cases}
      (x_1,\ldots,x_n,x_{n+1}^{\beta\alpha\beta^{-1}\alpha}) &\mbox{ if } x_1=\ldots=x_n=1 \\
      (x_1,\ldots,x_n,x_{n+1}) & \mbox{ if some } x_1,\ldots,x_{n-1} \neq 1,\,x_n=1 \\
      (x_1,\ldots,x_n,x_{n+1}) & \mbox{ if } x_1=\ldots=x_{n-1}=1,\, x_n\neq 1 \\
      (x_1,\ldots,x_n,x_{n+1}) & \mbox{ if some } x_1,\ldots,x_{n-1} \neq 1,\,x_n\neq 1 
    \end{cases}\\
 &= \begin{cases}
      (x_1,\ldots,x_n,x_{n+1}^{(1\ldots k)}) &\mbox{ if } x_1=\ldots=x_n=1 \\
      (x_1,\ldots,x_n,x_{n+1}) & \mbox{otherwise }  
    \end{cases}
\end{align*}
Thus we see that $\Sigma =  TG(n+1,(1,\ldots,k),1)$. 

We have shown that we can realise the two Toffoli gates that generate all of $\{TG(n+1,\alpha,o)\mid \alpha \in S_A\}$.
Thus by Corollary \ref{corrBnGen}, $B_{n+1}(A)$ is realised for all $n$, so by induction, all of $B(A)$ is realised.
\hfill$\Box$

Thus we have a small generating set for the bijections.

\begin{cor}
  Let $A =\{1,\dots,k\}$ be of odd order $k$.
  Let $\alpha = (1\ldots k), \beta=(1 2) \in S_A$.
 Then $B(A)$  is  
 realised by $\{\alpha,\beta,TG(2,\alpha,1), TG(2,\beta,1) \}$.
\end{cor}

We see that when $A$ is of odd order, $TG(i+1,\alpha,1)$ is in the reversible \toffalg{}  generated by $\{TG(i,\alpha,1)\mid \alpha \in S_A\}$,
which is not the case for even order $A$.
We see a distinct property looking at closure with temporary storage in the following.
Note that $TG(n,\alpha,o) = \mu((2,3,\ldots,n+1),\kappa(1,o,TG(n+1,\alpha,o)))$ so $TG(n,\alpha,o)$ 
is realised with strong temporary storage from $TG(n+1,\alpha,o)$.
We see that closure under $T_S$ is a stronger form of generation.


\begin{thm}[\cite{toff80} Thm 5.3]
\label{thmtoffsts}
Let $A$ be a set, $o\in A$. For every $n$, for every $\alpha$, $TG(n,\alpha,o)$ can be realised with strong temporary storage by 
 $\{TG(i,\alpha,o) \mid i \leq 3,\, \alpha \in S_A\}$.
\end{thm}
Proof:
We proceed by induction. Assume that we can construct   $TG(n-1,\alpha,o)$ for all $\alpha$.
Let  $p\in A$ some non-$o$ element, $\beta = (o\,p)$ be an transposition of $A$.
Define $f\in B_{n+1}(A)$ by
\begin{align*}
 f = ( TG(n-1,\beta,o) \oplus i_2) \bullet (i_{n-2} \oplus  TG(3,\alpha,o)    ) \bullet (TG(n-1,\beta,o) \oplus i_2)
\end{align*}
Then
\begin{align*}
 &f(x_1,\dots,x_{n+1}) \\
 &= \begin{cases}
      ( TG(n-1,\beta,o) \oplus i_2) \bullet (i_{n-2} \oplus  TG(3,\alpha,o)    )(x_1,\dots,x_{n-1}^\beta,x_n,x_{n+1})\\
	 \hspace{50mm}  \mbox{ if }x_1=\dots=x_{n-2}=o\\
      ( TG(n-1,\beta,o) \oplus i_2) \bullet (i_{n-2} \oplus  TG(3,\alpha,o)    )(x_1,\dots,x_{n-1},x_n,x_{n+1})\\
	\hspace{50mm} \mbox{otherwise}
    \end{cases}\\
 &= \begin{cases}
      ( TG(n-1,\beta,o) \oplus i_2)(x_1,\dots,x_{n-1}^\beta,x_n,x_{n+1}^\alpha)\\
	\hspace{30mm} \mbox{ if }x_1=\dots=x_{n-2}=x_{n-1}^\beta=x_n=o\\
      ( TG(n-1,\beta,o) \oplus i_2)(x_1,\dots,x_{n-1}^\beta,x_n,x_{n+1})\\
	\hspace{30mm}\mbox{ if } x_1=\dots=x_{n-2}=o \wedge (x_{n-1}^\beta \neq o \vee x_n \neq o)\\
      ( TG(n-1,\beta,o) \oplus i_2)(x_1,\dots,x_{n-1},x_n,x_{n+1}^\alpha)\\
	\hspace{30mm} \mbox{ if some }x_1,\dots,x_{n-2} \neq o,x_{n-1}=x_n=o\\
      ( TG(n-1,\beta,o) \oplus i_2)(x_1,\dots,x_{n-1},x_n,x_{n+1}) \hspace{10mm} \mbox{otherwise}
    \end{cases}
 \end{align*}
 \begin{align*}
 &= \begin{cases}
      (x_1,\dots,x_{n-1},x_n,x_{n+1}^\alpha) &\mbox{ if } x_1=\dots=x_{n-2}=x_{n-1}^\beta=x_n=o\\
      (x_1,\dots,x_{n-1},x_n,x_{n+1}) &\mbox{ if } x_1=\dots=x_{n-2}=o \,\wedge \\& \hspace{4mm}(x_{n-1}^\beta \neq o \vee x_n \neq o)\\
      (x_1,\dots,x_{n-1},x_n,x_{n+1}^\alpha) & \mbox{ if some }x_1,\dots,x_{n-2} \neq o,x_{n-1}=x_n=o\\
      (x_1,\dots,x_{n-1},x_n,x_{n+1})& \mbox{otherwise}
    \end{cases}
\end{align*}
Note that $f_i(x_1,\dots,x_{n+1}) =x_i$ for all $i\leq n$.

Let $g$ be the reduct 
\begin{align*}
 g(x_1,\dots,x_{n}) = \mu((1,2,\ldots,n-2,n,n+1),\kappa(n-1,p,f))
\end{align*}

From above we obtain, knowing that $p^\beta=o$ and $p\neq o$
\begin{align*}
 &\kappa(n-1,p,f) (x_1,\dots,x_{n})\\
 &= f(x_1,\ldots,x_{n-2},p,x_{n-1},x_n) \\
 &= \begin{cases}
      (x_1,\dots,x_{n-2},p,x_{n-1},x_{n}^\alpha) & \mbox{ if }x_1=\dots=x_{n-2}=p^\beta=x_{n-1}=o\\
      (x_1,\dots,x_{n-2},p,x_{n-1},x_{n}) &\mbox{ if } x_1=\dots=x_{n-2}=o \,\wedge\\ &\hspace{4mm}(p^\beta \neq o \vee x_{n-1} \neq o)\\
      (x_1,\dots,x_{n-2},p,x_{n-1},x_{n}^\alpha) & \mbox{ if }x_j\neq o \exists j\leq n-2,\\ &\hspace{4mm} p=x_{n-1}=o  \mbox{ [contradiction]}\\
      (x_1,\dots,x_{n-2},p,x_{n-1},x_{n}) & \mbox{otherwise}
    \end{cases} \\
 &= \begin{cases}
      (x_1,\dots,x_{n-2},p,x_{n-1},x_{n}^\alpha) & \mbox{ if }x_1=\dots=x_{n-2}=x_{n-1}=o\\
      (x_1,\dots,x_{n-2},p,x_{n-1},x_{n}) & \mbox{ if } x_1=\dots=x_{n-2}=o \wedge ( x_{n-1} \neq o)\\
      (x_1,\dots,x_{n-2},p,x_{n-1},x_{n}) & \mbox{otherwise}
    \end{cases}\\
 &= \begin{cases}
      (x_1,\dots,x_{n-2},p,x_{n-1},x_{n}^\alpha) & \mbox{ if }x_1=\dots=x_{n-2}=x_{n-1}=o\\
      (x_1,\dots,x_{n-2},p,x_{n-1},x_{n}) & \mbox{otherwise}
    \end{cases}
\end{align*}

Thus
\begin{align*}
 g(x_1,\dots,x_{n})
 &= \mu((1,2,\ldots,n-2,n,n+1),\kappa(n-1,p,f))(x_1,\dots,x_{n}) \\
 &= \begin{cases}
      (x_1,\dots,x_{n-2},x_{n-1},x_{n}^\alpha) & \mbox{ if }x_1=\dots=x_{n-2}=x_{n-1}=o\\
      (x_1,\dots,x_{n-2},x_{n-1},x_{n}) & \mbox{otherwise}
    \end{cases}\\
    &= TG(n,\alpha,o)
\end{align*}

Note that $f_{n-1}(x_1,\ldots,x_{n+1}) = x_{n-1}$, so we have strong temporary storage, thus completing our claim.
\hfill$\Box$

Thus we see some essential differences between isomorphic realisation and realisation with (strong) temporary storage.
In the next section we will further investigate the differences between certain types of realisation and closure.

\section{Various Closures}
\label{secclosures}

The following result shows us how the various class and closure operators 
that we have seen above relate and thus
we will be able to determine  information about the 
relationships between various types of closure.

\begin{thm}
 $KK=K$, $SS=S$, $CC=C$, $C_\Delta C_\Delta=C_\Delta$, $SC=CS$, $KC=CK$, 
 $C_\Delta K=KC_\Delta$, $C_\Delta S=SC_\Delta$ and $SK=KS$.
\end{thm}
Proof:
The $K$ and $S$ operators are idempotent as a simple implication of the definitions.
The $C$ and $C_\Delta$ operators are idempotent because they are algebraic closure operations.

To show that $SC=CS$ we show that all \toffalg{} operations commute with $s$. 
We use the operation set $\{\oplus,\pi_\alpha,\circ_k\mid \alpha \in S_{(\N)},\, k\in \N\}$.
 We start with the inclusion $SC \subseteq CS$.
 Let $f\in M_{n,m}(A)$, $g\in M_{l,p}(A)$, $r\in \N$, $I \in \{1,\dots,m+p\}^r$ without repetitions, $\alpha \in S_m$.

Let $J^\prime\subseteq \{1,\dots,r \}$ such that $j\in J^\prime$ iff $I_j \leq m$.
Write $J^\prime=\{j_1,\dots,j_t\}$ in ascending order $j_1<\dots<j_t$, then define $J=(j_1,\dots,j_t) \in \{1,\dots,r\}^t$. 
Then let $I^\prime=(I_{j_i}\mid i \in \{1,\dots,t\}) \in \{1,\dots,m\}^t$.
Let $\{1,\dots,r\} \setminus J = \{\bar j_1,\dots,\bar j_u\}$. 
Then let $I^{\prime\prime} = (I_{\bar j_i}-m \mid i \in \{1,\dots,u\}) \in \{1,\dots,p\}^u$.
Finally define $\beta \in S_r$ with 
\[ \beta: i \mapsto h \mbox{ such that } \begin{cases}
                     h=j_i & \mbox{ if }i\leq t \\
                     h = \bar j_{i-t}  &\mbox{ if } i > t
                    \end{cases}
\] which can be written as
\[\beta = \pmtv[1\ldots {t}{t+1}\ldots {r}] {{1}{j_1} {t}{j_t} {t+1}{\bar j_1} {r}{\bar j_u} }\]

Note that $i\leq t$ iff $I_{\beta(i)} \leq m$ and that $I_i = I^\prime_{\beta^{-1}i}$ if $I_i\leq m$,
$I_i=I^{\prime\prime}_{(\beta^{-1}i)-t}+m$ otherwise.
Then we claim that
\begin{align*}
 \mu(I,f\oplus g) &= \pi_\beta (\mu(I^\prime,f) \oplus \mu(I^{\prime\prime},g)) 
 \end{align*}
The left hand side is
\begin{align*}
 \mu(I,f\oplus g)_i &=
 \begin{cases}
  (f)_{I_i} &\mbox{ if } I_i \leq m \\
  (g)_{I_i-m} & \mbox{ if }I_i > m
 \end{cases}
\end{align*}
The right hand side is
\begin{align*}
 \pi_\beta (\mu(I^\prime,f) \oplus \mu(I^{\prime\prime},g))_i 
 &= (\mu(I^\prime,f) \oplus \mu(I^{\prime\prime},g))_{\beta^{-1}i} \\
 &=\begin{cases}
    \mu(I^\prime,f)_{\beta^{-1}i} & \mbox{ if }{\beta^{-1}i}\leq t \\
    \mu(I^{\prime\prime},g)_{(\beta^{-1}i)-t} & \mbox{ if }{\beta^{-1}i} > t
   \end{cases}\\
 &=\begin{cases}
    (f)_{I^\prime_{\beta^{-1}i}} & \mbox{ if }{\beta^{-1}i}\leq t \\
    (g)_{I^{\prime\prime}_{(\beta^{-1}i)-t}} & \mbox{ if }{\beta^{-1}i} > t
   \end{cases}\\
 &=\begin{cases}
    (f)_{I_{i}} & \mbox{ if }I_i\leq m \\
    (g)_{I_{i}-m} &\mbox{ if } I_{i} > m
   \end{cases}
\end{align*}
which shows our claim.

 It is a simple calculation that
 \begin{align*}
 \mu(I,\pi_\alpha f) &= \mu(\alpha^{-1}(I),f)
 \end{align*}
 where $\alpha^{-1}$ acts upon the entries in $I$, so $(\alpha^{-1}I)_i = \alpha^{-1} (I_i)$.

For composition, we use a similar argument to the $\oplus$ case above. 
Assume that $I$ is increasing. 
Let $t$ be such that $I_t < m$ and $I_{t+1} > m$ and let 
 $I^\prime = (I_1,\dots,I_t)$.
Let $I^{\prime\prime} = (1,\dots,k)\oplus(I_{t+1}-m+k,\dots,I_r-m+k)$.
We claim that
 \begin{align*}
 \mu(I,f\circ_k g) &= \mu(I^\prime,f) \circ_k \mu(I^{\prime\prime},g)
\end{align*}
The left hand side is
\begin{align*}
 \mu(I,f\circ_k g)_i &=
 \begin{cases}
  f_{I_i} \circ_k g &\mbox{ if } i \leq t\\
  g_{I_i-m+k} & \mbox{ if }i > t
 \end{cases}
\end{align*}
The right hand side is 
\begin{align*}
 (\mu(I^\prime,f) \circ_k \mu(I^{\prime\prime},g))_i &= 
 \begin{cases}
  \mu(I^\prime,f)_i \circ_k \mu(I^{\prime\prime},g) &\mbox{ if } i \leq t\\
  \mu(I^{\prime\prime},g)_{i-t+k}	& \mbox{ if }i > t
 \end{cases}\\
 &= 
 \begin{cases}
  f_{I^\prime_i} \circ_k g & \mbox{ if }i \leq t \;\;\mbox{because of the first $k$ entries in }I^{\prime\prime}\\
  g_{I^{\prime\prime}_{i-t+k}}	& \mbox{ if }i > t
 \end{cases}\\
 &= 
 \begin{cases}
  f_{I_i} \circ_k g &\mbox{ if } i \leq t \\
  g_{I_{i}-m+k}	&\mbox{ if } i > t
 \end{cases}
\end{align*}
Which is what we wanted. 
If $I$ is not increasing, then there is a permutation $\beta \in S_r$ such that $\pi_\beta I$ is increasing.
Then 
\[\mu(I,f \circ_k g) = \pi_{\beta^{-1}} \mu(\pi_\beta I,f \circ_k g) 
= \pi_{\beta^{-1}}( \mu(I^\prime,f) \circ_k \mu(I^{\prime\prime},g))
\]
Thus we see that $SC \subseteq CS$.

For the converse, i.e.\ $CS \subseteq SC$, we use similar techniques.
Let $f\in M_{n,m}(A)$, $g\in M_{l,p}(A)$, $t,u\in \N$, $I^\prime \in \{1,\dots,m\}^t$, $I^{\prime\prime} \in \{1,\dots,p\}^u$, 
$\alpha \in S_t$.

Let $I=I^\prime \oplus (I^{\prime\prime}+m)$, i.e.\ $I_i = I^\prime_i$ for $i\leq t$, $I_i = I^{\prime\prime}_{i-t}+m$ for $i>t$.
Then it is a simple calculation to see that
\begin{align*}
 \mu(I^\prime,f) \oplus \mu(I^{\prime\prime},g) = \mu(I,f\oplus g)
\end{align*}

Let $I \in \{1,\dots,m\}^t$, define $I^\alpha_i = I_{\alpha^{-1}i}$.
Then we see that
\begin{align*}
 \pi_\alpha \mu(I,f) = \mu(I^\alpha,f)
\end{align*}

Let $\beta \in S_p$ be defined by $\beta^{-1}(i) = I^{\prime\prime}_i$ for $i\leq u$, 
with the rest filled in to make it a permutation. Because $I^{\prime\prime}$ contains no repeats,
we know this can be done.
Let $k \leq \min(n,u)$.
Then we claim that
\begin{align*}
 \mu(I^\prime,f) \circ_k \mu(I^{\prime\prime},g) = \mu(I^\prime \oplus (m+1,\dots,m+u-k),f \circ_k (\pi_\beta \circ_p g))
\end{align*}

For $i \leq t$ we have
\begin{align*}
 &(\mu(I^\prime,f) \circ_k \mu(I^{\prime\prime},g))_i(x_1,\dots,x_{l+n-k}) \\
   &= f_{I^\prime_i} \circ_k \mu(I^{\prime\prime},g) (x_1,\dots,x_{l+n-k})\\
   &= f_{I^\prime_i}(g_{I^{\prime\prime}_1}(x_1,\dots,x_{l}),\dots,g_{I^{\prime\prime}_k}(x_1,\dots,x_{l}),x_{l+1},\dots,x_{l+n-k}) \\
   &= f_{I^\prime_i}(g_{\beta^{-1}1}(x_1,\dots,x_{l}),\dots,g_{\beta^{-1}k}(x_1,\dots,x_{l}),x_{l+1},\dots,x_{l+n-k}) \\
   &= (\mu(I^\prime \oplus (m+1,\dots,m+u-k),f \circ_k (\pi_\beta \circ_p g)))_i(x_1,\dots,x_{l+n-k})
\end{align*}
while for $t < i \leq m+u-k$ we have
\begin{align*}
 &(\mu(I^\prime,f) \circ_k \mu(I^{\prime\prime},g))_i(x_1,\dots,x_{l+n-k}) \\
 &= g_{I^{\prime\prime}_{k+i-t}} (x_1,\dots,x_l) \\
 &= g_{\beta^ {-1}(k+i-t)} (x_1,\dots,x_l) \\
 &= (\pi_\beta \circ_p g)_{k+(i-t)}(x_1,\dots,x_l) \\
 &= (f \circ_k (\pi_\beta \circ_p g))_{m+i-t}(x_1,\dots,x_{l+n-k}) \\
 &= \mu(I^\prime \oplus (m+1,\dots,m+u-k),f \circ_k (\pi_\beta \circ_p g))_i(x_1,\dots,x_{l+n-k})
\end{align*}

Thus we have that $CS \subseteq SC$ and thus $SC=CS$.

Similarly we show that $KC=CK$.
We start with $KC \subseteq CK$.
Let $f\in M_{n,m}(A)$, $g\in M_{l,p}(A)$, $i \in \{1,\dots,n+l\}$, $a\in A$. 
It is a simple application of the definitions to see that
\begin{align*}
 \kappa(i,a,f \oplus g) &= \begin{cases}
                           \kappa(i,a,f) \oplus g &\mbox{ if } i \leq n\\
                           f \oplus \kappa(i-n,a,g) & \mbox{ if }i > n
                          \end{cases}
 \end{align*}
Let $f\in M_{n,m}(A)$, $\alpha \in S_m$, $i \in \{1,\dots,n\}$, $a\in A$. 
It is a simple calculation that
\begin{align*}
 \kappa(i,a,\pi_\alpha f) &= \pi_\alpha \kappa(i,a,f)   
\end{align*}
Let $f\in M_{n,m}(A)$, $g\in M_{l,p}(A)$, $k \leq \min(n,p)$, $i \in \{1,\dots,n+l\}$, $a\in A$. 
Applying the definitions gives us
\begin{align*}
 \kappa(i,a,f \circ_k g) &= \begin{cases}
                           f \circ_k  \kappa(i,a,g) & \mbox{ if }i \leq n\\
                           \kappa(i+k-n,a,f) \circ_k g & \mbox{ if }i > n
                          \end{cases}
\end{align*}

For $CK \subseteq KC$ we proceed as follows.
Note that we must also include the nonapplication of $\kappa$ as a case here.
Let $f\in M_{n,m}(A)$, $g\in M_{l,p}(A)$, $i_1 \in \{1,\dots,n\}$, $i_2\in \{1,\dots,l\}$, $a_1,a_2\in A$. 
Then we claim that
\begin{align*}
 \kappa(i_1,a_1,f) \oplus g &= \kappa(i_1,a_1,f\oplus g)\\
 f \oplus \kappa(i_2,a_2,g) &= \kappa(i_2+n,a_2,f\oplus g)\\
 \kappa(i_1,a_1,f) \oplus \kappa(i_2,a_2,g) &= \kappa(i_1,a_1,\kappa(i_2+n,a_2,f\oplus g))
\end{align*}
The first two follow by simple application of the definitions.
The third claim combines these two.

As we saw above, $\pi_\alpha$ commutes with $\kappa$.

Let $f\in M_{n,m}(A)$, $g\in M_{l,p}(A)$, $k \leq \min(n,p)$, $i_1 \in \{1,\dots,n\}$, $i_2\in \{1,\dots,l\}$, $a_1,a_2\in A$. 
Let $\beta \in S_n$ be $(i_1,i_1+1,\dots,n)$. Then we claim that
\begin{align*}
f \circ_k \kappa(i_2,a_2,g) &= \kappa(i_2,a_2,f\circ_k g) \\
 \kappa(i_1,a_1,f) \circ_k g &= 
 \begin{cases}
  \kappa(l+n-k,a_1,(f \circ_n \pi_\beta) \circ_k g) &\mbox{ if } i_1 \leq k \\
  \kappa(i_1+l-k,a_1,f \circ_k g) & \mbox{ if }i_1 > k
 \end{cases}\\
 \kappa(i_1,a_1,f) \circ_k \kappa(i_2,a_2,g) &= 
 \begin{cases}
  \kappa(i_2,a_2,\kappa(l+n-k,a_1,(f \circ_n \pi_\beta) \circ_k g)) &\mbox{ if } i_1 \leq k \\
  \kappa(i_2,a_2,\kappa(i_1+l-k,a_1,f \circ_k g)) &\mbox{ if } i_1 > k
 \end{cases}
\end{align*}
The first claim is direct from the definitions.
The second requires some more work. Let $i\in \{1,\dots,m+p-k\}$.
\begin{align*}
 &(\kappa(i_1,a_1,f) \circ_k g )_i(x_1,\ldots,x_{n+l-k-1})\\
 &= \begin{cases}
     (\kappa(i_1,a_1,f)_i \circ_k g)(x_1,\ldots,x_{n+l-k-1}) & \mbox{ if }i\leq m\\
     g_{i-m+k}(x_1,\dots,x_l) &\mbox{ if } i > m
    \end{cases}\\
 &= \begin{cases}
     f_i(g_1(x_1,\dots,x_l),\dots,g_{i_1-1}(\dots),a_1,g_{i_1}(\dots),\dots\\
       \hspace{18mm}\dots,g_k(x_1,\dots,x_l),x_{l+1},\dots,x_{n+l-k-1}) & \mbox{ if }i_1 \leq k,\, i\leq m\\
     f_i(g_1(x_1,\dots,x_l),\dots,g_k(x_1,\dots,x_l),x_{l+1},\dots\\
        \hspace{18mm}\dots,x_{i_1-k+l-1},a_1,x_{i_1+l-k},\dots,x_{n+l-k-1}) &\mbox{ if } i_1 > k,\, i\leq m\\
     g_{i-m+k}(x_1,\dots,x_l) & i > m
    \end{cases}\\
 &= \begin{cases}
     (f\circ_n\pi_\beta)_i(g_1(x_1,\dots,x_l),\dots,g_{i_1-1}(\dots),g_{i_1}(\dots),\dots\\
       \hspace{18mm}\dots,g_k(x_1,\dots,x_l),x_{l+1},\dots,x_{n+l-k-1},a_1) &\mbox{ if }i_1 \leq k,\, i\leq m\\
     \kappa(i_1+l-k,a_1,f\circ_k g)_i(x_1,\ldots,x_{n+l-k-1}) & \mbox{ if }i_1 > k,\, i\leq m\\
     g_{i-m+k}(x_1,\dots,x_l) & \mbox{ if }i > m
    \end{cases}\\
 &= \begin{cases}
     \kappa(l+n-k,a_1,(f \circ_n \pi_\beta)\circ_k g)_i(x_1,\ldots,x_{n+l-k-1}) &\mbox{ if }i_1 \leq k,\, i\leq m\\
     \kappa(i_1+l-k,a_1,f\circ_k g)_i(x_1,\ldots,x_{n+l-k-1}) &\mbox{ if } i_1 > k,\, i\leq m\\
     g_{i-m+k}(x_1,\dots,x_l) &\mbox{ if } i > m
    \end{cases}\\
 &= \begin{cases}
     \kappa(l+n-k,a_1,(f \circ_n \pi_\beta)\circ_k g)_i(x_1,\ldots,x_{n+l-k-1}) &\mbox{ if } i_1 \leq k\\
     \kappa(i_1+l-k,a_1,f\circ_k g)_i(x_1,\ldots,x_{n+l-k-1}) &\mbox{ if } i_1 > k
    \end{cases}
\end{align*}
which is what we wanted.
The third case is the combination of the first two cases.

So every expression in $CK$ can be written in $KC$, so we obtain $CK=KC$.

We see that $SK=KS$ because it makes no difference whether the inputs are fed constants and then some outputs are ignored, 
or some outputs are ignored and then some inputs are fed constants.
To see this formally,
let $f\in M_{n,m}(A)$, $1\leq i\leq n$, $r\in \N$, $I\in \{1,\ldots,m\}^r$  and $a\in A$.
\begin{align*}
 \mu(I,\kappa(i,a,f))(x_1,\ldots,x_{n-1}) &= \mu(I,f(x_1,\ldots,x_{i-1},a,x_i,\ldots,x_{n-1})) \\
   &= (f_j(x_1,\ldots,x_{i-1},a,x_i,\ldots,x_{n-1}) \mid j\in I) \\
   &= (\kappa(i,a,f_j)\mid j\in I)(x_1,\ldots,x_{n-1}) \\
   &= \kappa(i,a,\mu(I,f))(x_1,\ldots,x_{n-1})
\end{align*}

Lastly, we show that $C_\Delta$ is well behaved by showing that $\Delta$ and $\nabla$ commute with $S$ and $K$.
Let $f\in M_{n,m}(A)$, $r\in \N$, $I\in \{1,\ldots,n\}^r$ and $a\in A$.
\begin{align*}
 \mu(I,\Delta f)(x_1,\dots,x_{n-1}) &= \mu(I,f(x_1,x_1,x_2,\ldots,x_{n-1}))\\
  &= (f_j(x_1,x_1,x_2,\ldots,x_{n-1}) \mid j\in I)\\
  &= \Delta (f_j(x_1,x_2,\ldots,x_{n-1}) \mid j\in I)\\
  &= \Delta \mu(I,f)(x_1,\dots,x_{n-1})\\
 \mu(I,\nabla f)(x_1,\dots,x_{n+1}) &= \mu(I, f(x_2,\dots,x_{n+1})) \\
  &= (f_j(x_2,\dots,x_{n+1})\mid j\in I) \\
  &= \nabla \mu(I,f)(x_1,\dots,x_{n+1})
\end{align*}     
So we see that $C_\Delta S = S C_\Delta$. 

Now let $f\in M_{n,m}(A)$, $1\leq i\leq n$ and $a\in A$.
\begin{align*}
&\kappa(i,a,\Delta f)(x_1,\dots,x_{n-2})\\ &=  
               (\Delta f)(x_1,\dots,x_{i-1},a,x_i,\dots,x_{n-2}) \\
  &= f(x_1,x_1,x_2,\dots,x_{i-1},a,x_i,\dots,x_{n-2}) \\
  &= \begin{cases}
      f(a,a,x_1,\ldots,x_{n-2}) &\mbox{ if }i=1 \\
      f(x_1,x_1,a,x_2,\ldots,x_{n-2}) &\mbox{ if }i=2 \\
      f(x_1,x_1,x_2,\ldots,x_{i-1},a,x_i,\ldots x_{n-2}) &\mbox{ if } i>2
     \end{cases}\\
  &= \begin{cases}
      \kappa(1,a,\kappa(1,a,f)) & \mbox{ if }i=1 \\
      \Delta \kappa(i+1,a,f) &\mbox{ if } i\geq2 
     \end{cases}\\
\end {align*}
\begin{align*}
 \kappa(i,a,\nabla f)(x_1,\dots,x_n) &= (\nabla f)(x_1,\dots,x_{i-1},a,x_i,\dots,x_n) \\
  &= 
  \begin{cases}
   (\nabla f)(a,x_1,\dots,x_n)	&\mbox{ if } i=1\\
   (\nabla f)(x_1,x_2,\dots,x_{i-1},a,x_i,\dots,x_n) &\mbox{ if } i>1
  \end{cases}\\
  &= 
  \begin{cases}
   f(x_1,\dots,x_n)	& \mbox{ if }i=1\\
   f(x_2,\dots,x_{i-1},a,x_i,\dots,x_n) & \mbox{ if }i>1
  \end{cases}\\
  &= 
  \begin{cases}
   f(x_1,\dots,x_n)	&\mbox{ if } i=1\\
   \nabla \kappa(i-1,a,f)(x_1,\dots,x_n) &\mbox{ if } i>1
  \end{cases} 
\end{align*}     
Thus we have $KC_\Delta \subseteq C_\Delta K$. 
For the converse we  proceed as follows.
Suppose $n=1$, then $\Delta \kappa(1,a,f) = f(a) \in KC_\Delta$ and we are done.
Suppose $n=2$, then $\kappa(1,a,f)$ and $\kappa(2,a,f)$ are both unary and thus fixed by $\Delta$ so we are done.
Suppose $n\geq 3$, we calculate
\begin{align*}
 &\hspace{-10mm}\kappa(i,a,f) (x_1,\ldots,x_{n-1})\\
  &= \begin{cases}
       f(a,x_1,x_2,\ldots,x_{n-1}) & \mbox{ if }i=1 \\
       f(x_1,a,x_2,\ldots,x_{n-1}) &\mbox{ if } i=2 \\
       f(x_1,x_2,\ldots,x_{i-1},a,x_i,\ldots,x_{n-1}) &\mbox{ if } i>2
     \end{cases} \\
 &\hspace{-10mm}\Rightarrow \Delta \kappa(i,a,f) (x_1,\ldots,x_{n-2})\\
   &= \begin{cases}
      f(a,x_1,x_1,x_2\ldots,x_{n-2}) & \mbox{ if }i=1 \\
       f(x_1,a,x_1,x_2,\ldots,x_{n-2}) &\mbox{ if } i=2 \\
       f(x_1,x_1,x_2,\ldots,x_{i-2},a,x_{i-1},\ldots,x_{n-2}) & \mbox{ if }i>2
     \end{cases} \\
  &= \begin{cases}
      \kappa(2,a,\Delta(f\bullet \pi_{(1\,2\,3)}))(x_1,\dots,x_{n-2}) & \mbox{ if }i=1 \\
      \kappa(2,a,\Delta (f\bullet \pi_{(2\,3)}))(x_1,\dots,x_{n-2})  &\mbox{ if } i=2 \\
      \kappa(i-1,a,\Delta f)(x_1,\dots,x_{n-2}) &\mbox{ if }i>2
     \end{cases} \\
&\hspace{-10mm}\nabla \kappa(i,a,f)(x_1,\dots,x_n)\\ &= \kappa(i,a,f)(x_2,\dots,x_n)\\
  &= f(x_2,\dots,x_i,a,x_{i+1},\dots,x_n)\\
  &= (\nabla f)(x_1,\dots,x_i,a,x_{i+1},\dots,x_n) \\
  &= \kappa(i+1,a,\nabla f)(x_1,\dots,x_n)
\end{align*}
So we see that 
$K C_\Delta \subseteq C_\Delta K \subseteq K C_\Delta$ so they are equal.
\hfill$\Box$

\begin{figure}
\begin{center}
\begin{tikzpicture}[scale=.7]
\node (20) at (2,4) {$KSC(F)$};
\node (21) at (-2,1)  {$T(F)$};
\node (22) at (0,-1)  {$T_S(F)$};
\node (10) at (2,-3) {$C(F)$};
\node (KSC) at (3,-1) {$\cdot$};
\node (KC) at (2,1) {$KC(F)$};
\node (SC) at (4,1) {$SC(F)$};
\node (RKSC) at (-3,2) {$RKSC(F)$};
\node (RT) at (-4,-1) {$RT(F)$};
\node (RTS) at (-2,-3) {$RT_S(F)$};
\node (RC) at (0,-5) {$RC(F)$};

\node (KSCD) at (5,6) {$KSC_\Delta(F)$};
\node (KCD) at  (5,3) {$KC_\Delta(F)$};
\node (SCD) at  (8,3) {$SC_\Delta(F)$};
\node (KcSCD) at (8,0.5) {$\cdot$};
\node (CD) at   (7,-1) {$C_\Delta(F)$};

\draw (20) -- (21) -- (22) -- (10);
\draw (20) -- (KC) -- (KSC);
\draw (20) -- (SC) -- (KSC) -- (10)--(RC);
\draw (20)--(RKSC)--(RT)--(RTS)--(RC);
\draw (21)--(RT);
\draw (22)--(RTS);

\draw (KSCD)--(20);
\draw (KSCD)--(KCD)--(KcSCD);
\draw (KSCD)--(SCD)-- (KcSCD)--(CD)--(10);
\draw (KCD)--(KC);
\draw (SCD)--(SC);
\end{tikzpicture}
\begin{tikzpicture}[scale=.7]
\node (20) at (2,4) {$KSC(F)$};
\node (22) at (-1,-2)  {$T_S(F)=RT_S(F)$};
\node (10) at (2,-4) {$C(F)=RC(F)$};
\node (KC) at (2,0) {$KC(F)$};
\node (SC) at (5,0) {$SC(F)$};
\node (RKSC) at (-1,2) {$RKSC(F)$};
\node (RT) at (-1.5,0) {$T(F)=RT(F)$};

\draw (22) -- (10);
\draw (20) -- (KC) -- (10);
\draw (20) -- (SC) -- (10);
\draw (20)--(RKSC)--(RT)--(22);

\end{tikzpicture}

\end{center}
\caption{At the top we have inclusion of the various closure operations applied to a set of maps $F \subseteq M(A)$. 
When we are looking at a set $F \subseteq B(A)$ on a finite $A$ and are not interested in $\Delta$, we get the inclusions in the bottom diagram. 
Note
that we know that $SC(F)\cap KC(F)=C(F)$ in this case.}
\label{figureinclusion}
\end{figure}

Let $F\subseteq M(A)$ be a collection of mappings.
Then $g \in SKC(F)$ is equivalent to saying that $g$ is realised by $F$.
We see that we have some inclusions amongst the various closure operations introduced above, 
obtaining the inclusion diagram in Figure \ref{figureinclusion}.

We know that many of these inclusions are strict.
For $A$ of even order, we know from Corollary \ref{thmtoffeven} and Theorem \ref{thmtoffsts} that $C(F)$ is strictly included in $T_S(F)$.

On the other hand, for specific classes of $F$, some closure classes fall together.
If $F \subseteq B(A)$, then $RC(F)=C(F)$, $RT_S(F) = T_S(F)$ and $RT(F)=T(F)$. 
Thus we obtain the second inclusion diagram, where we also omit the $\Delta$ operator.

One of the general problems is to look at the ways in which we can 
define these various classes as a form of closure via a Galois connection.
This has been preliminarily investigated in \cite{jerebek} for $C(F)$ and $T(F)$ from reversible $F$.

\section{Conclusion}

We have presented a language based upon clone theory in order to discuss systems of 
mappings $A^m\rightarrow A^n$ on a set $A$.
We have then written Toffoli's model of reversible computation as a theory of mapping composition, 
finding that larger alphabets leads to a more complex system than binary reversible logic.

We see that the ideas can be used more generally. We are now confronted with 
the spectrum of questions that arise naturally
in clone theory and related fields of general algebra, such as the size and 
structure of the lattice of multiclones and reversible \toffalgs{} on a given set.
It would be of value to determine a set of exioms for a \toffalg{} or a reversible \toffalg, 
this might imprive the clarity of many of the proofs.
The largest challenge is to determine a suitable combinatorial 
structure to be used for invariance type results of the Pol-Inv type for the 
five natural types of closure.
It seems probable that the results in \cite{aaronson,jerebek} will be of use to develop such a theory.

\section{Acknowledgments}
I would like to thank Erhard Aichinger for many interesting conversations and suggestions about developing
and presenting this paper.
Erkko Lehtonen made many excellent and detailed suggestions that have vastly improved the paper
and for which I am extremely grateful.

\bibliography{reversible}

\def\cprime{$'$}
\begin{thebibliography}{10}

\bibitem{aaronson}
S.~Aaronson, D.~Grier, and L.~Schaeffer.
\newblock (2015).
\newblock The classification of reversible bit operations.
\newblock {\em Electronic Colloquium on Computational Complexity}, (66).

\bibitem{stefanescu13}
J.A. Bergstra, C.A. Middleburg, and G.~{\c S}tef{\u a}nescu.
\newblock (2013).
\newblock Network algebra for synchronous dataflow.
\newblock https://arxiv.org/pdf/1303.0382.pdf.

\bibitem{boykettkariville17}
T.~Boykett, J.~Kari, and V.~Salo.
\newblock (2017).
\newblock Finite generating sets for reversible gate sets under general
  conservation laws.
\newblock {\em Theoretical Computer Science}, pages~--.
\newblock Electronic publication Feb 2017.

\bibitem{BS}
S.~Burris and H.~P. Sankappanavar.
\newblock (1981).
\newblock {\em A course in universal algebra}, volume~78 of {\em Graduate Texts
  in Mathematics}.
\newblock Springer-Verlag, New York-Berlin.

\bibitem{coucierolehtonen}
M.~Couceiro and E.~Lehtonen.
\newblock (2012).
\newblock Galois theory for sets of operations closed under permutation,
  cylindrification, and composition.
\newblock {\em Algebra Universalis}, 67(3):273--297.

\bibitem{GAP4}
The GAP~Group.
\newblock (2013).
\newblock {\em {GAP -- Groups, Algorithms, and Programming, Version 4.6.5}}.

\bibitem{engineering}
N.~Gershenfeld.
\newblock (1996).
\newblock Signal entropy and the thermodynamics of computation.
\newblock {\em IBM Systems Journal}, 35(3 \& 4).

\bibitem{jensen14}
R.~Jensen.
\newblock (2014).
\newblock The logic of reversible structures.
\newblock Master's thesis, University of Copenhagen.

\bibitem{jerebek}
E.~Je{\v{r}}{\'a}bek, (September 2014).
\newblock Answer to classifying reversible gates.
\newblock Theoretical Computer Science Stack Exchange.
\newblock http://cstheory.stackexchange.com/questions/25730, Accessed June
  2015.

\bibitem{lehtonen}
E.~Lehtonen.
\newblock (2010).
\newblock Closed classes of functions, generalized constraints, and clusters.
\newblock {\em Algebra Universalis}, 63(2-3):203--234.

\bibitem{malcev}
A.~I. Mal{\cprime}tsev.
\newblock (1976).
\newblock {\em Iterativnye algebry {P}osta}.
\newblock Novosibirsk. Gos. Univ., Novosibirsk.
\newblock Seriya ``Biblioteka Kafedry Algebry i Matematicheskoi Logiki
  Novosibirskogo Universiteta'', Vyp. 16. [``Library of the Department of
  Algebra and Mathematical Logic of Novosibirsk State University'' Series, No.
  16].

\bibitem{physics}
M.~A. Nielsen and I.~L. Chuang.
\newblock (2000).
\newblock {\em Quantum computation and quantum information}.
\newblock Cambridge University Press, Cambridge.

\bibitem{pk79}
R.~P{\"o}schel and L.~A. Kalu{\v{z}}nin.
\newblock (1979).
\newblock {\em Funktionen- und {R}elationenalgebren [Function and Relation
  algebras]}, volume~15 of {\em Mathematische Monographien [Mathematical
  Monographs]}.
\newblock VEB Deutscher Verlag der Wissenschaften, Berlin.
\newblock Ein Kapitel der diskreten Mathematik. [A chapter in discrete
  mathematics].

\bibitem{stefanescu2000}
G.~{\c S}tef{\u a}nescu.
\newblock (2000).
\newblock {\em Network algebra}.
\newblock Discrete Mathematics and Theoretical Computer Science (London).
  Springer-Verlag London, Ltd., London.

\bibitem{szendrei}
{\'A}.~Szendrei.
\newblock (1986).
\newblock {\em Clones in universal algebra}, volume~99 of {\em S\'eminaire de
  Math\'ematiques Sup\'erieures [Seminar on Higher Mathematics]}.
\newblock Presses de l'Universit\'e de Montr\'eal, Montreal, QC.

\bibitem{toff80}
T.~Toffoli.
\newblock (1980).
\newblock Reversible computing.
\newblock Technical Report MIT/LCS/TM-151, MIT.

\bibitem{tofflncs}
T.~Toffoli.
\newblock (1980).
\newblock Reversible computing.
\newblock In {\em Automata, languages and programming ({P}roc. {S}eventh
  {I}nternat. {C}olloq., {N}oordwijkerhout, 1980)}, volume~85 of {\em Lecture
  Notes in Comput. Sci.}, pages 632--644. Springer, Berlin-New York.

\bibitem{yangetal}
G.~Yang, X.~Song, M.~Perkowski, and J.~Wu.
\newblock (2005).
\newblock Realizing ternary quantum switching networks without ancilla bits.
\newblock {\em J. Phys. A}, 38(44):9689--9697.

\end{thebibliography}

\end{document}